\documentclass[journal,twoside,web]{ieeecolor}

\usepackage{generic}

\usepackage{amsmath,amssymb,amsfonts,mathtools,amsthm}
\usepackage{cite}
\usepackage{algorithmic}
\usepackage{graphicx}
\usepackage{algorithm,algorithmic}
\usepackage{hyperref}
\hypersetup{hidelinks=true}
\usepackage{textcomp}
\def\BibTeX{{\rm B\kern-.05em{\sc i\kern-.025em b}\kern-.08em
    T\kern-.1667em\lower.7ex\hbox{E}\kern-.125emX}}
\newtheorem{theorem}{Theorem}[section]
\newtheorem{proposition}[theorem]{Proposition}
\newtheorem{lemma}[theorem]{Lemma}
\theoremstyle{definition}

\newtheorem{remark}[theorem]{Remark}
\newtheorem{assumption}[theorem]{Assumption}
\newtheorem{example}[theorem]{Example}
\newtheorem{corollary}[theorem]{Corollary}

\newcommand{\R}{\mathbb{R}}

\newcommand{\argmin}{\mathop{\mathrm{arg\,min}}}

\newcommand{\rvW}{\mathsf{w}}
\newcommand{\rvX}{\mathsf{x}}
\newcommand{\rvU}{\mathsf{u}}
\newcommand{\rvY}{\mathsf{y}}
\newcommand{\rvZ}{\mathsf{z}}
\newcommand{\rvV}{\mathsf{v}}
\newcommand{\rvN}{\mathsf{n}}
\newcommand{\rvC}{\mathsf{c}}

\newcommand{\mrvW}{\mathsf{W}}
\newcommand{\mrvX}{\mathsf{X}}
\newcommand{\mrvU}{\mathsf{U}}

\newcommand{\rcd}{P}

\begin{document}

\title{Risk-Sensitive Stochastic Optimal Control on Positive Cones}

\author{
Mohamad Al Ahdab and John Leth
\thanks{Mohamad Al Ahdab and John-Josef Leth are with the
Nonlinear and Optimal Control Lab, Section of Automation and Control,
Department of Electronic Systems, Aalborg University,
Aalborg, Denmark
(e-mail: \{maah, jjl\}@es.aau.dk).}
}

\maketitle

\begin{abstract}
We study risk-sensitive stochastic optimal control problems in which the state is constrained to a proper cone. Specifically, we use the entropic-risk criterion to cover risk-averse, risk-seeking, and risk-neutral (expected-cost) problems. We focus on systems for which the Bellman operator maps linear functions to affine functions of the state and show how this structure reduces the optimal control problem to a finite-dimensional recursion on the dual cone. We derive optimal value functions and feedback policies for finite horizons and characterize the optimal time-average cost and stationary policy through the smallest fixed point of the recursion. Additionally, under a finite-step cost observability assumption, we establish uniqueness of the fixed point, geometric convergence of the recursion, and finite-horizon performance bounds for risk-seeking and risk-neutral problems. We also establish uniform-in-time bounds on the expected state norm for risk-averse and risk-neutral problems. Moreover, we show that risk-averse problems with disturbances entering the dynamics and costs affinely admit an equivalent deterministic infsup formulation. We provide several examples satisfying the affine Bellman property, including stochastic positive systems and systems on the positive-semidefinite cone, and show how Linear Quadratic (Exponential) Gaussian problems are related to the latter.

\end{abstract}


\section{Introduction}

Finding explicit solutions to optimal control problems remains challenging in general. The linear-quadratic regulator (LQR) is one of the best-known classes of problems for which an explicit optimal value function and feedback policy can be obtained. In addition to LQR, the structural properties of positive and monotone systems have been exploited to simplify the solution of various optimal control problems \cite{Colaneri_Middleton_Blanchini_2016,Dhingra_Colombino_Jovanovic_Rantzer_Smith_2016,Zorzan_Rantzer_2017}. More recently, explicit solutions and finite-dimensional characterizations of linear optimal value functions have been established for classes of positive linear systems and related problems \cite{Rantzer_2022,DynamicProgramming_YuchaoRantzer,Ohlin_Pates_Arcak_2025,semilineardynamicprogramminganalysis}. Although LQR and optimal control for positive systems appear to be different problems, both can be related using linear value functions on positive cones \cite{Pates_Rantzer_2024}. Specifically, the quadratic value function of LQR becomes linear after lifting the state to the cone of positive-semidefinite (PSD) matrices. This observation motivates the formulation of optimal control problems on general proper cones in \cite{Pates_Rantzer_2024}, which brings these classes of problems within a common framework and may allow other problems admitting linear optimal value functions to be identified.

In many applications, however, the dynamics are affected by disturbances and uncertainty. This raises the question of whether the linear value-function structure can also be exploited to obtain explicit solutions in such settings. For LQR, the $H_\infty$ (minimax) formulation provides an approach for accounting for worst-case disturbances \cite{Glover_Doyle_1988,Basar_Bernhard_2008}. Similarly, minimax optimal control problems for positive systems with linear costs have been studied in \cite{Gurpegui_Tegling_Rantzer_2023,Gurpegui_Jeeninga_Tegling_Rantzer_2026}. Related minimax adaptive control problems on cones have also been considered in \cite{Rantzer_2025_AdaptiveCones}. Alternatively, when a probabilistic description of the disturbances is available, one can formulate a stochastic optimal control problem. In this setting, different realizations of the disturbances may lead to different state trajectories and incurred costs. The resulting cost is thus a random variable, and a criterion for evaluating its distribution is needed to compare different control policies. Typically, this is done by considering the expected value of the cost, as in the classical linear-quadratic-Gaussian (LQG) problem \cite{Whittle_2002}. For stochastic positive systems, \cite{Colaneri_Middleton_Blanchini_2016} establishes convexity of the expected terminal-state cost for positive Markovian bilinear systems. Furthermore, \cite[Secs.~5 and~6]{semilineardynamicprogramminganalysis} establishes linear optimal value functions and certainty-equivalence properties for positive semilinear systems with multiplicative uncertainty (discounted and undiscounted total costs) and Markov jump parameters (undiscounted total costs). Explicit solutions for continuous-time positive systems with multiplicative noise are derived in \cite{SDEAlba} for finite-horizon and discounted infinite-horizon expected-cost problems.

For systems operating over long periods of time in the presence of persistent disturbances, the undiscounted total cost may be infinite even under a stabilizing control policy. While discounting future costs can address this issue, it assigns progressively smaller weights to costs incurred further into the future. An alternative is to consider the time-average (also referred to as the long-run average or ergodic) cost criterion, which evaluates the asymptotic cost per unit time without discounting future costs. This criterion is relevant for long-term operation in applications such as resource management and process control. Although the aforementioned works establish favorable structural properties and explicit solutions for several classes of stochastic positive systems, they do not directly address the time-average optimal control problem within the general proper-cone framework of \cite{Pates_Rantzer_2024}.

Furthermore, minimizing the expected value of a random cost does not distinguish between cost distributions having the same mean. An alternative is to use a \emph{risk-sensitive} criterion, such as the entropic-risk criterion \cite{Howard_Matheson_1972,Whittle_1990,Whittle_2002}. Specifically, for a nonnegative random cost $\rvC$, this criterion takes the form $\frac{1}{\theta}\log\mathbb{E}[\exp(\theta\rvC)]$, where $\theta\neq0$ is a risk parameter describing the preference toward uncertainty. For $\theta>0$, minimizing the criterion is equivalent to minimizing $\mathbb{E}[\exp(\theta\rvC)]$, in which the penalty $\exp(\theta\rvC)$ grows exponentially as the realized cost increases. Thus, $\theta>0$ provides a \emph{risk-averse} perspective. In contrast, for $\theta<0$, minimizing the criterion instead maximizes $\mathbb{E}[\exp(\theta\rvC)]$, in which $\exp(\theta\rvC)$ becomes a reward that increases exponentially as the realized cost decreases and approaches zero as the realized cost increases. Therefore, the reward gained from a possible very low cost can outweigh the reward lost from a possible high cost, and hence $\theta<0$ provides a \emph{risk-seeking} perspective \cite{Whittle_2002}.\footnote{The sign convention of the risk parameter in our paper is opposite to the sign convention adopted in Whittle \cite[p.~6]{Whittle_2002}.} For example, consider an input $u\in\{0,1\}$ chosen before observing $\rvZ$, with cost $\rvC(u)=5u+\rvZ(1-u)$, where $\rvZ=10$ with probability $1/2$ and $\rvZ=0$ otherwise. The risk-averse decision is $u=1$, since it incurs a guaranteed cost of $5$ and avoids the possibility of paying $10$. On the other hand, the risk-seeking decision is $u=0$, since it favors the possibility of paying $0$ despite the risk of paying $10$. Note that both inputs have the same expected cost of $5$, so if the criterion is the expectation (\emph{risk-neutral}) instead, then both decisions are optimal. Moreover, if $\mathbb{E}[\exp(\varepsilon\rvC)]<\infty$ for some $\varepsilon>0$, then, for sufficiently small $|\theta|$, the criterion is approximated by $\mathbb{E}[\rvC]+\frac{\theta}{2}\operatorname{Var}(\rvC)$ \cite[p.~7]{Whittle_2002}. 

Entropic-risk optimal control has been studied for finite-state Markov decision processes, linear-exponential-quadratic-Gaussian (LEQG) control (LQG with the entropic risk criterion), and more general stochastic systems \cite{Howard_Matheson_1972,Jacobson_1973,Whittle_1981,Whittle_1990,Bauerle_Rieder_2014}. Finite-horizon risk-averse problems on Borel spaces are studied in \cite{Chapman_RADP2022}, while time-average entropic-risk optimal control on general state spaces, including problems with unbounded nonnegative costs, is considered in \cite{Jaskiewicz_2007}. Entropic-risk optimal control also has established connections with minimax control. In particular, connections between controllers satisfying an $H_\infty$-norm bound and LEQG were investigated in~\cite{Glover_Doyle_1988}. Related connections with deterministic dynamic games have been studied in~\cite{Jacobson_1973,Campi_James_1996,Basar_Bernhard_2008}. More generally, the Donsker-Varadhan variational formula provides a connection between the entropic-risk criterion and distributionally robust formulations, where alternative probability distributions are considered and deviations from the nominal distribution are penalized by their relative entropy~\cite{Dupuis_James_Petersen_2000,Hansen_Sargent_2008}.

Motivated by the deterministic formulation in \cite{Pates_Rantzer_2024} and the above developments in stochastic and risk-sensitive optimal control, we consider discrete-time stochastic optimal control on proper cones under the entropic-risk criterion. Specifically, we focus on problems whose Bellman operator maps linear functions to affine functions of the state. The consideration of affine rather than linear functions is motivated by the fact that disturbances may give rise to a nonzero cost even at the origin. Our main contrinutions are as follows
\begin{enumerate}
\item[(i)] We characterize finite-horizon optimal value functions and policies through a finite-dimensional recursion on the dual cone (Proposition~\ref{prop:finite_horizon_affine}) and establish properties of its effective domain (Proposition~\ref{prop:rs_domain}).

\item[(ii)] We establish equivalent conditions for fixed-point existence and characterize the optimal time-average cost and stationary policy through the smallest fixed point (Theorem~\ref{thm:average_cost}). We also provide an optimization characterization and total-cost results (Corollary~\ref{cor:supersolution_program} and Proposition~\ref{prop:stochastic_total_cost}).

\item[(iii)] Under finite-step cost observability, we establish fixed-point uniqueness, geometric convergence, a convex optimization characterization, and a uniform bound on the finite-horizon suboptimality gap of the stationary optimal policy for nonpositive risk parameters (Theorem~\ref{thm:detectability_consequences}).

\item[(iv)] Under finite-step cost observability, we establish uniform-in-time bounds on the expected state norm for nonnegative risk parameters and exponential-moment bounds for positive risk parameters (Proposition~\ref{prop:finite_step_state_bounds}).

\item[(v)] For risk-averse problems with dynamics and costs affine in the disturbances, we establish an equivalent deterministic infsup formulation (Section~\ref{sec:infsup_connection}).

\item[(vi)] We provide examples satisfying the affine Bellman property and recover LQG and LEQG from the PSD formulation which extends the deterministic construction in \cite{Pates_Rantzer_2024} (Section~\ref{sec:PSDexample}).
\end{enumerate}

\section{Notation and Preliminaries}
\label{sec:Notation}

Euclidean spaces are equipped with the Euclidean norm $\|\cdot\|$ and their Borel $\sigma$-algebras, and $\operatorname{int}\mathcal S$ denotes the interior of a set $\mathcal S$.
We write $\mathbb R_{\geq0}^d$ and $\mathbb R_{>0}^d$ for the vectors in $\mathbb R^d$ with nonnegative and positive entries, and $\mathbb S^n$, $\mathbb S_{\geq}^n$, and $\mathbb S_{>}^n$ for the symmetric, PSD, and positive-definite $n\times n$ matrices. For vectors and matrices, $|\cdot|$ and inequalities are understood entrywise, and $a^{(i)}$ and $A^{(ij)}$ denote the $i$th entry of vector $a$ and the $(i,j)$th entry of matrix $A$.

We fix a probability space $(\Omega,\mathcal F,\mathbb P)$, and equalities and inequalities involving random variables hold $\mathbb P$-almost surely unless stated otherwise. Random variables are denoted by sans-serif symbols, such as $\rvW$, and their realizations by the corresponding lowercase letters, such as $w$.

A Markov kernel from a measurable space $\mathcal Z$ to $\mathcal W\subseteq\mathbb R^{n_w}$ is a map $\rcd:\mathcal B(\mathcal W)\times\mathcal Z\to[0,1]$ such that $\rcd(\cdot\mid z)$ is a probability measure on $\mathcal W$ for each $z\in\mathcal Z$, and $z\mapsto \rcd(B\mid z)$ is measurable for each Borel set $B\in\mathcal B(\mathcal W)$. Moreover, $\rcd$ is a regular conditional distribution of $\rvW$ given $\rvZ$ if $\rcd(B\mid\rvZ)=\mathbb P(\rvW\in B\mid\rvZ)$ for every Borel set $B$, where $\rcd(B\mid \rvZ)$ is considered the map $z\mapsto \rcd(B\mid \rvZ=z)$. In this case, we write $\rvW\mid\rvZ=z\sim \rcd(\cdot\mid z)$, or briefly $\rvW\mid z\sim \rcd(\cdot\mid z)$, and $\mathbb E[\varphi(\rvW)\mid z]:=\int_{\mathcal W}\varphi(w)\,\rcd(dw\mid z)$ for integrable $\varphi$.

Fix $\theta\in\mathbb R$. For a real-valued random variable $\rvY$ and a sub-$\sigma$-algebra $\mathcal G\subseteq\mathcal F$, define the entropic-risk functional by
\begin{equation}
\label{eq:risk_functional}
\mathcal R_\theta(\rvY\mid\mathcal G):=
\begin{cases}
\mathbb E[\rvY\mid\mathcal G], & \theta=0,\\
\dfrac{1}{\theta}\log\mathbb E[\exp(\theta\rvY)\mid\mathcal G], & \theta\neq0,
\end{cases}
\end{equation}
whenever the corresponding quantities are well defined, and write $\mathcal R_\theta(\rvY)$ when $\mathcal G=\{\varnothing,\Omega\}$. We also use \eqref{eq:risk_functional} for nonnegative random variables that may take the value $+\infty$, with the conventions $e^{\theta(+\infty)}=+\infty$ for $\theta>0$, $e^{\theta(+\infty)}=0$ for $\theta<0$, and $\log0=-\infty$. Conditioning on a random variable $\rvZ$ means conditioning on $\sigma(\rvZ)$, and $\mathcal R_\theta(\rvY\mid z)$ and $\mathbb E[\rvY\mid z]$ denote the values at $\rvZ=z$, also when conditioning on several random variables, as in $\mathcal R_\theta(\cdot\mid x,u)$.

The functional $\mathcal R_\theta$ is monotone,
\begin{equation}
\label{eq:monoticity_of_R}
\rvY_1\leq\rvY_2\quad\Rightarrow\quad\mathcal R_\theta(\rvY_1\mid \mathcal{G})\leq\mathcal R_\theta(\rvY_2\mid \mathcal{G}),
\end{equation}
and translation invariant, i.e., $\mathcal R_\theta(\rvY+\rvC\mid\mathcal G)=\mathcal R_\theta(\rvY\mid\mathcal G)+\rvC$ for every $\mathcal G$-measurable $\rvC$. Moreover, the tower property of conditional expectation gives, for $\mathcal G\subseteq\mathcal H\subseteq\mathcal F$ and $\mathcal H$-measurable $\rvY_0$,
\begin{equation}
\label{eq:recursive_R}
\mathcal R_\theta\left(\rvY_0+\mathcal R_\theta(\rvY_1\mid\mathcal H)\mid\mathcal G\right)=\mathcal R_\theta\left(\rvY_0+\rvY_1\mid\mathcal G\right),
\end{equation}
whenever the corresponding quantities are well defined.

A set $\mathcal K\subseteq\mathbb R^d$ is a cone if $\alpha x\in\mathcal K$ for every $x\in\mathcal K$ and $\alpha\geq0$. Its dual cone is $\mathcal K^\ast:=\{\lambda\in\mathbb R^d:\lambda^\top x\geq0\text{ for all }x\in\mathcal K\}$. A cone $\mathcal K$ is proper if it is closed, convex, pointed, i.e., $\mathcal K\cap(-\mathcal K)=\{0\}$, and has nonempty interior. In this case, $\mathcal K^\ast$ is also proper. We write $y\geq_{\mathcal K}x$ if $y-x\in\mathcal K$. For any nonempty cone $\mathcal K$, we have
\begin{equation}
\label{eq:rep}
\gamma+\beta^\top x\geq0~\text{for all }x\in\mathcal K
\Longleftrightarrow~\gamma\geq0\text{ and }\beta\in\mathcal K^\ast,
\end{equation}
for every $\gamma\in\mathbb R$ and $\beta\in\mathbb R^d$.

For a proper cone $\mathcal K$, we use the following properties:
\begin{subequations}
\label{eq:cone_interior_properties}
\begin{equation}
\label{eq:cone_span_bidual}
\mathcal K^{\ast\ast}=\mathcal K,\qquad\operatorname{span}\mathcal K=\mathbb R^d.
\end{equation}
Moreover,
\begin{equation}
\label{eq:cone_interior_addition}
\alpha\operatorname{int}\mathcal K+\mathcal K\subseteq\operatorname{int}\mathcal K,~\alpha>0.
\end{equation}
For every $p\in\operatorname{int}\mathcal K$ and $v\in\mathbb R^d$, there exist $\varepsilon>0$ and $\mu\geq1$ such that
\begin{equation}
\label{eq:cone_interior_comparison}
p-tv\in\mathcal K,~0\leq t\leq\varepsilon,\qquad\mu p-v\in\mathcal K.
\end{equation}
Finally, for every $\lambda\in\operatorname{int}\mathcal K^\ast$, there exists $\delta>0$ such that
\begin{equation}
\label{eq:cone_interior_norm}
\lambda^\top x\geq\delta\|x\|,~x\in\mathcal K.
\end{equation}
\end{subequations}
The properties in \eqref{eq:cone_interior_addition} and \eqref{eq:cone_interior_comparison} follow from the definition of interior and the cone property \cite[Sec.~2.4.1]{BoydVandenberghe2004}. To obtain \eqref{eq:cone_interior_norm}, choose $\delta>0$ such that $\lambda-\delta v\in\mathcal K^\ast$ for every $\|v\|\leq1$, and then take $v=x/\|x\|$ for $x\neq0$. These properties also hold in finite-dimensional real inner-product spaces.

\section{Problem Formulation and Solutions}
In this section, we formulate the stochastic optimal control problem and introduce the assumptions defining the class of interest. We then characterize the solutions to the finite-horizon and time-average infinite-horizon problems and address the infinite-horizon total-cost problem.

\subsection{Problem Setup and Assumptions}
\label{sec:ProbAss}
Let $\mathcal X\subseteq\mathbb R^{n_x}$ be a proper cone, with dual cone $\mathcal X^\ast$, which is also proper. Given a Borel measurable map $f:\mathbb R^{n_x}\times\mathbb R^{n_u}\times\mathbb R^{n_w}\to\mathbb R^{n_x}$, a Borel measurable stage cost $\ell:\mathbb R^{n_x}\times\mathbb R^{n_u}\times\mathbb R^{n_w}\to\mathbb R$, a set-valued map $\mathcal U:\mathcal X\to2^{\mathbb R^{n_u}}$ with a Borel graph, a Borel set $\mathcal W\subseteq\mathbb R^{n_w}$, an initial condition $x_0\in\mathcal X$, a vector $\lambda\in\mathcal{X}^\ast$, and a sequence (policy) $\pi_N:=\{a_t\}_{t=0}^{N-1}$, $N\geq1$ of Borel measurable feedback decision rules $a_t:\mathcal X\to\mathbb R^{n_u}$, define
\begin{equation}
\label{eq:Nstage_cost_noterminal}
J_{N,\lambda}^{\pi_N}(x_0):=\mathcal R_\theta\left(S_{N,\pi_N}^\lambda(x_0)\right),
\end{equation}
where
\begin{equation}
\label{eq:Nstage_cost_terminal}
S_{N,\pi_N}^\lambda(x_0):=\displaystyle\sum_{t=0}^{N-1}\ell(\rvX_t,\rvU_t,\rvW_t)+\lambda^\top\rvX_N,
\end{equation}
with $\rvX_t,~\rvU_t$, and $\rvW_t$ defined through
\begin{subequations}
\label{eq:general_dynamics}
\begin{align}
\label{eq:policy}
&\rvU_t=a_t(\rvX_t),\\
\label{eq:disturbance}
&\begin{multlined}
\rvW_t\mid(\rvX_t,\rvU_t)=(x,u)\sim \rcd(\cdot\mid x,u),
\end{multlined}\\
\label{eq:dynamics}
&\rvX_{t+1}=f(\rvX_t,\rvU_t,\rvW_t),\quad \rvX_0=x_0
\end{align}
\end{subequations}
for $t\in\{0,\ldots,N-1\}$, and  where
\begin{multline}
\rcd:\mathcal B(\mathcal W)\times\left\{(x,u)\in\mathcal X\times\mathbb R^{n_u}:u\in\mathcal U(x)\right\}\to[0,1]
\end{multline}
is a regular conditional distribution of $\rvW_t$ given
$(\rvX_t,\rvU_t)$. Furthermore, given $(\rvX_t,\rvU_t)$, the disturbance $\rvW_t$ is conditionally independent of all past states, controls, and disturbances.
More precisely, for every $B\in\mathcal B(\mathcal W)$,
$\mathbb P\bigl(\rvW_t\in B\mid(\rvX_s,\rvU_s)_{s=0}^{t},(\rvW_s)_{s=0}^{t-1}\bigr)=\rcd(B\mid\rvX_t,\rvU_t)$.
The expectation in $\mathcal R_\theta$ is taken with respect to the probability law induced by
$x_0$, $\pi_N$, and $\rcd$. For $N=0$, the sum in \eqref{eq:Nstage_cost_terminal} is empty by definition. Additionally, we will write $J_N^{\pi_N}(x_0):=J_{N,0}^{\pi_N}(x_0)$ in the rest of the paper for ease of notation. Moreover, we write $\pi:=\{a_t\}_{t=0}^{\infty}$ for infinite-horizon policies.

For each $x\in\mathcal X$, define the admissible input set
\begin{multline*}
\mathcal A(x):=\Bigl\{u\in\mathcal U(x):\\\rcd\bigl(\{w\in\mathcal W:f(x,u,w)\in\mathcal X\}\mid x,u\bigr)=1\Bigr\}.
\end{multline*}
A finite- or infinite-horizon policy is admissible if $a_t(x)\in\mathcal A(x)$ for every $x\in\mathcal X$ and $t\geq 0$. Let $\Pi_N$ and $\Pi$ denote the sets of admissible $N$-stage and infinite-horizon policies, respectively. For $\pi\in\Pi$, we use $S_{N,\pi}^\lambda$, $J_{N,\lambda}^\pi$, and $J_N^\pi$ for the corresponding quantities associated with its first $N$ decision rules.

We seek to characterize the optimal values and policies of the finite-horizon and time-average infinite-horizon problems, respectively:
\begin{subequations}
\label{eq:control_problems}
\begin{align}
&\inf_{\pi_N\in\Pi_N}J_{N,\lambda}^{\pi_N}(x_0),\label{eq:finite_horizon_ocp}\\
J_{\mathrm{av}}^\star(x_0)&:=\inf_{\pi\in\Pi}\limsup_{N\to\infty}\frac{1}{N}J_N^\pi(x_0).\label{eq:average_cost_ocp}
\end{align}
\end{subequations}

We now specify the class of problems considered in this paper. Unless stated otherwise, the assumptions introduced in this subsection are taken to hold from the point at which they are introduced. Additional assumptions are imposed explicitly where needed.
\begin{assumption}
\label{ass:basic}
For every $x\in\mathcal X$, the admissible input set $\mathcal A(x)$ is nonempty. Moreover, for every $x\in\mathcal X$ and $u\in\mathcal A(x)$,
$\ell(x,u,\rvW)\geq0$.
\end{assumption}

Under Assumption~\ref{ass:basic}, the accumulated costs in \eqref{eq:Nstage_cost_terminal} are nonnegative for admissible policies. We allow $J_{N,\lambda}^{\pi_N}(x_0)=+\infty$ when the expected cost is infinite for $\theta=0$, or when the exponential moment is infinite for $\theta>0$. Such policies remain admissible.

The problem class are
expressed through the one-stage Bellman operator, $T_\theta$, which is defined as
\begin{equation}
    \label{eq:Bellman_operator}
    (T_\theta J)(x) := \!\!\!\!
    \inf_{u\in\mathcal A(x)}\!\!\!\! \mathcal{R}_\theta
    (\mathcal{H}_J(x,u,\rvW)\mid x,u)
\end{equation}
with $$\mathcal{H}_J(x,u,w):= \ell(x,u,w)+J(f(x,u,w))$$
and $J: \mathcal{X}\to \mathbb{R}\cup \{\infty\}$. Note that $T_\theta$ is monotone in the following sense:
\begin{multline}\label{eq:mono}
\text{
If $J(x)\leq \tilde{J}(x)$ for all $x\in\mathcal{X}$, then}\\ \text{$(T_\theta J)(x)\leq (T_\theta \tilde J )(x)$ for all $x\in\mathcal{X}$.}
\end{multline}
Indeed, for any $x\in\mathcal{X}$ and $u\in\mathcal{A}(x)$, we have $f(x,u,\rvW)\in\mathcal{X}$, hence $J(f(x,u,\rvW))\leq \tilde J(f(x,u,\rvW))$ by assumption. Adding the same stage cost $\ell$, applying $\mathcal{R}_\theta$ and using \eqref{eq:monoticity_of_R}, and minimizing over the same admissible set $\mathcal{A}(x)$ then gives \eqref{eq:mono}.
Observe also that the operator $T_\theta$ is translation invariant, i.e., for a constant $C(x)\equiv c\in\mathbb{R}$, we have $(T_\theta(J+C))(x)=(T_\theta J)(x)+c$.
Now starting from some $J_{0,\lambda}(x)=\lambda^\top x$ for $\lambda\in\mathcal{X}^{\ast}$, define for $N\geq 1$ the following
\begin{multline}
\label{eq:rs_JN1}
J_{N,\lambda}(x):=(T_\theta J_{N-1,\lambda})(x)\\
=\!\!\!\inf_{u\in\mathcal A(x)}\!\!\!\!
\mathcal{R}_\theta
(\mathcal{H}_{J_{N-1,\lambda}}(x,u,\rvW)\mid x,u).
\end{multline}
The following lemma relates \eqref{eq:rs_JN1} to the finite-horizon problem in \eqref{eq:finite_horizon_ocp} \footnote{see \cite{Bertsekas_Shreve_1996,Bauerle_Rieder_2014} for general dynamic-programming background and \cite{Chapman_RADP2022} for a related bounded-cost entropic-risk formulation on Borel spaces}.
\begin{lemma}
\label{lem:finite_horizon_verification}
Let $N\geq1$ and $\lambda\in\mathcal X^\ast$. Suppose that $J_{k,\lambda}$ in \eqref{eq:rs_JN1} is finite, nonnegative, and Borel measurable for $k\in\{0,\ldots,N\}$. Then
$J_{N,\lambda}(x)\leq J_{N,\lambda}^{\pi_N}(x),~x\in\mathcal X,~\pi_N\in\Pi_N$.
Moreover, the identity
\begin{align}
\label{eq:rs_JN2}
J_{N,\lambda}(x)=\inf_{\pi_N\in\Pi_N} J_{N,\lambda}^{\pi_N}(x)
\end{align}
holds whenever $\pi_N=\pi_N^\star:=\{a_t^\star\}_{t=0}^{N-1}\in\Pi_N$ with
\begin{equation}
\label{eq:rs_optimal_feedback}
a_t^\star(x)\in\argmin_{u\in\mathcal A(x)}\mathcal R_\theta\left(\mathcal H_{J_{N-t-1,\lambda}}(x,u,\rvW)\mid x,u\right),
\end{equation}
for every $x\in\mathcal X$ and
$t\in\{0,\ldots,N-1\}$.
\end{lemma}
The proof is given in Appendix~\ref{sec:finite_horizon_verification_proof}. Henceforth, we will write $J_N=J_{N,0}$ for ease of notation.

For $\theta=0$, we impose the following
integrability assumption.
\begin{assumption}
    \label{ass:finite_expectation}
    For every $x\in\mathcal{X}$, $u\in\mathcal{A}(x)$, and $\lambda\in\mathcal{X}^\ast$, we have
        \begin{equation}  \mathcal{R}_0\left(\mathcal{H}(x,u,\rvW,\lambda)\mid x,u\right)<\infty
    \end{equation}
with $\mathcal{H}(x,u,w,\lambda):=\ell(x,u,w)+\lambda^\top f(x,u,w)$ (that is, $\mathcal H=\mathcal H_J$ with $J(x)=\lambda^Tx$).
\end{assumption}
Note that $\mathcal{H}(x,u,\rvW,\lambda)\geq 0$ every $x\in\mathcal{X}$, $u\in\mathcal{A}(x)$, and $\lambda\in\mathcal{X}^\ast$ by assumption \ref{ass:basic}. It follows that 
\begin{align}\label{eq:Rtheta_geq0}
\mathcal{R}_\theta(\mathcal{H}(x,u,\rvW,\lambda)\mid x,u)\geq 0,
\end{align}
which for $\theta\geq 0$ is clear, and for $\theta<0$ we have $0<e^{\theta\mathcal{H}(x,u,\rvW,\lambda)}\leq1$ so the corresponding logarithmic moment is nonnegative after division by $\theta$. Next, we state the following convexity/concavity property of $\lambda\mapsto\mathcal{R}_\theta(\mathcal{H}(x,u,\rvW,\lambda\mid x,u)$. 
\begin{proposition}
\label{prop:convexthetaR}
Let $x\in\mathcal X$, $u\in\mathcal A(x)$, $\lambda,\widetilde\lambda\in\mathcal X^\ast$, $\tau\in[0,1]$, and define $\lambda_\tau:=\tau\lambda+(1-\tau)\widetilde\lambda$. Then,\footnote{with the convention $0\cdot\infty=0$.} for $\theta\neq 0$,
\begin{multline}
\label{eq:convexthetaR}
\theta\mathcal R_\theta(\mathcal H(x,u,\rvW,\lambda_\tau)\mid x,u)\leq
\tau\theta\mathcal R_\theta(\mathcal H(x,u,\rvW,\lambda)\mid x,u)\\+(1-\tau)\theta\mathcal R_\theta(\mathcal H(x,u,\rvW,\widetilde\lambda)\mid x,u).
\end{multline}
Additionally, for $\theta=0$, we have
\begin{multline}
\label{eq:expectationR}
\mathcal R_0(\mathcal H(x,u,\rvW,\lambda_\tau)\mid x,u)=
\tau\mathcal R_0(\mathcal H(x,u,\rvW,\lambda)\mid x,u)\\+(1-\tau)\mathcal R_0(\mathcal H(x,u,\rvW,\widetilde\lambda)\mid x,u),
\end{multline}
\end{proposition}
\begin{proof}
For $\theta=0$, \eqref{eq:expectationR} follows from the linearity of expectation. Let $\theta\neq 0$. For $\tau\in\{0,1\}$, we have $\lambda_\tau\in\{\lambda,\widetilde\lambda\}$, and \eqref{eq:convexthetaR} holds with equality. For $\tau\in(0,1)$, let $\rvY:=\mathcal H(x,u,\rvW,\lambda)$ and $\widetilde\rvY:=\mathcal H(x,u,\rvW,\widetilde\lambda)$. Since $\lambda\mapsto\mathcal H(x,u,\rvW,\lambda)=\ell(x,u,\rvW)+\lambda^\top f(x,u,\rvW)$ is affine, we have $\mathcal H(x,u,\rvW,\lambda_\tau)=\tau\rvY+(1-\tau)\widetilde\rvY$. Applying H\"older's inequality, with $p=1/\tau$ and $q=1/(1-\tau)$, to the positive factors $(e^{\theta\rvY})^{\tau}$ and $(e^{\theta\widetilde\rvY})^{1-\tau}$ gives
\begin{multline*}
\mathbb E\left[\exp(\theta\mathcal H(x,u,\rvW,\lambda_\tau))\mid x,u\right]\leq\\
\left(\mathbb E\left[\exp(\theta\rvY)\mid x,u\right]\right)^{\tau}\left(\mathbb E\left[\exp(\theta\widetilde\rvY)\mid x,u\right]\right)^{1-\tau}.
\end{multline*}
Taking logarithms, and using $\theta\mathcal R_\theta(\cdot\mid x,u)=\log\mathbb E[\exp(\theta\,\cdot)\mid x,u]$ from \eqref{eq:risk_functional}, proves \eqref{eq:convexthetaR}.
\end{proof}

For the case when $\theta>0$, we have that $\mathcal{R}_\theta$ may not be finite for many common distributions (e.g., for a random variable having a probability distribution with exponantially decaying tails). Therefore, taking a finiteness assumption for $\theta>0$, similar to the assumption of the one in \eqref{ass:finite_expectation} for $\theta=0$, is too restrictive. To adress this, we define the following effective domain for $\lambda$ given a specific $\theta\in\mathbb{R}$ 
\begin{multline}
\label{eq:rs_effective_domain}
\mathcal D_\theta:=\Big\{
\lambda\in\mathcal X^\ast:\\
\inf_{u\in\mathcal A(x)}\!\!\!\!
\mathcal{R}_\theta(\mathcal{H} (x,u,\rvW,\lambda)
\mid x,u)
<\infty ~\forall x\in\mathcal X\Big\}.
\end{multline}
For the given assumptions so far, the following proposition provide the properties of $\mathcal D_{\theta}$ that will be essential for the subsequent derivations in the paper. 
\begin{proposition}
\label{prop:rs_domain} The following properties hold.
\begin{enumerate}
\item[(i)] If $\theta\leq0$, then $\mathcal D_\theta=\mathcal X^\ast$.
\item[(ii)] $\mathcal D_\theta$ is downward-closed with respect to the order induced on $\mathcal X^\ast$, i.e., if
$\lambda\in\mathcal D_\theta$, $\widetilde\lambda\in\mathcal X^\ast$, and
$\lambda\geq_{\mathcal X^\ast}\widetilde\lambda$, then
$\widetilde\lambda\in\mathcal D_\theta$.
\item[(iii)] Consequently, whenever $\mathcal D_\theta$ is nonempty, it contains
$0$ and is star-shaped about $0$, i.e., if $\lambda\in\mathcal D_\theta$ then $\tau\lambda\in\mathcal{D}_\theta$ for $\tau\in[0,1]$.
\end{enumerate}
\end{proposition}
\begin{proof}
Let $\lambda\in\mathcal X^\ast$, $x\in\mathcal X$, and $u\in\mathcal A(x)$. For $\theta=0$, $\mathcal R_\theta(\mathcal H(x,u,\rvW,\lambda)\mid x,u)$ is finite by Assumption~\ref{ass:finite_expectation}. For $\theta<0$, it is finite since $0<e^{\theta\mathcal{H}(x,u,\rvW,\lambda)}\leq1$, by $\mathcal H(x,u,\rvW,\lambda)\geq0$. Since $\mathcal A(x)$ is nonempty by Assumption~\ref{ass:basic}, the infimum in \eqref{eq:rs_effective_domain} is finite for every $x$, which proves (i).

To prove (ii), let $J(x)=\lambda^\top x$ and $\widetilde J(x)=\widetilde\lambda^\top x$. By definition of $\mathcal{D}_\theta$, $(T_\theta J)(x)<\infty$ for every $x\in\mathcal{X}$, and since $\widetilde J\leq J$ on $\mathcal X$, the monotonicity of $T_\theta$ in \eqref{eq:mono} gives $(T_\theta\widetilde J)(x)\leq(T_\theta J)(x)<\infty$, i.e., $\widetilde\lambda\in\mathcal{D}_\theta$. Finally, (iii) follows from (ii) with $\widetilde\lambda=0$ and with $\widetilde\lambda=\tau\lambda$, $\tau\in[0,1]$, since $\lambda=\tau\lambda+(1-\tau)\lambda\geq_{\mathcal X^\ast}\tau\lambda$.
\end{proof}

We now state the final assumption for the class of problems that this paper is tackling.  
\begin{assumption}
\label{ass:affine_bellman}
For a fixed $\theta\in\mathbb{R}$, and for every $\lambda\in\mathcal{D}_\theta$, there exist $\Gamma_\theta(\lambda)\in\mathbb R$, and $\Phi_\theta(\lambda)\in\mathbb R^{n_x}$, such that for $J(x)=\lambda^\top x$, we have
\begin{equation}
\label{eq:bellman_affinity}
\Gamma_\theta(\lambda)+\Phi_\theta(\lambda)^\top x \!=\!(T_{\theta}J)(x)\!=\!\!\!\inf_{u\in\mathcal A(x)}\!\!\! \mathcal{R}_\theta(\mathcal{H} (x,u,\rvW,\lambda)\mid x,u),
\end{equation}
for every \(x\in\mathcal X\). The minimum is attained by a stationary feedback decision rule
\begin{align}\label{eq:opt_sol}
x\longmapsto a(x;\lambda)\in\mathcal A(x).
\end{align}
where the occurrence of $\lambda$ in \eqref{eq:opt_sol} is to indicate that the optimal decision rule is parameterized by $\lambda$.
\end{assumption}
This assumption restrict the attention to a problem instance for which the one-stage Bellman minimization in \eqref{eq:bellman_affinity} preserves affine functions on $\mathcal X$.

We discuss several examples of classes of problems which satisfy \ref{ass:affine_bellman}. 

First, the following example considers a stochastic version of the positive linear systems with linear costs and state-dependent input constraints studied in \cite{Rantzer_2022,DynamicProgramming_YuchaoRantzer} with minimax formulations in \cite{Gurpegui_Tegling_Rantzer_2023,Gurpegui_Jeeninga_Tegling_Rantzer_2026}, and a risk-neutral continuous-time stochastic formulation with multiplicative noise in \cite{SDEAlba}. This class of problems have several applications, including voltage control in direct-current power networks \cite{Gurpegui_Tegling_Rantzer_2023}, water-network management \cite{Gurpegui_Jeeninga_Tegling_Rantzer_2026}, and portfolio allocation \cite{SDEAlba}.
\begin{example}
\label{ex:polyhedral_ex1intro}
Let $\mathcal X=\mathbb R_{\geq0}^{n_x}$ ($\mathcal X^\ast=\mathbb R_{\geq0}^{n_x}$) and $\mathcal U(x)=\{u\in\mathbb R^{n_u}:|u|\leq Ex\}$, where $E\geq0$, and consider
\begin{equation}
\label{eq:polyhedral_assumption_model}
f(x,u,w)=Ax+Bu+Fw,~\ell(x,u,w)=s^\top x+r^\top u,
\end{equation}
where $F\geq0$, $A-|B|E\geq0$, and $s\geq E^\top|r|$ (such that $s^\top x+r^\top u\geq0$ for all $x\in\mathcal{X}$ and $u\in\mathcal{U}(x)$). Consider $\rvW\in\mathbb R_{\geq0}^{n_d}$ with conditionally independent components and conditional mean $\mu_w(x,u)=\bar w+A_wx+B_wu>0$ for every $x\in\mathcal X$ and $u\in\mathcal U(x)$, where $\bar w>0$. We consider either the inverse-Gaussian (IG) model
\begin{equation}
\label{eq:polyhedral_ig}
\rvW^{(i)}\mid(x,u)\sim\operatorname{IG}\left(\mu_w^{(i)}(x,u),2\kappa^{(i)}(\mu_w^{(i)}(x,u))^2\right),
\end{equation}
where the parameters of $\operatorname{IG}$ are the mean and shape, respectively, and $\kappa^{(i)}>0$, or the Poisson model
\begin{equation}
\label{eq:polyhedral_poisson}
\rvW^{(i)}\mid(x,u)\sim\operatorname{Poisson}\left(\mu_w^{(i)}(x,u)\right),
\end{equation}
for $i\in\{1,\ldots,n_d\}$.

In both cases, $\mathcal{A}(x)=\mathcal{U}(x)$ for every $x\in\mathcal{X}$. Additionally, for $\lambda\in\mathcal{D}_\theta$, we have
\begin{equation}
\label{eq:polyhedral_ig_risk}
\mathcal R_\theta\left(\lambda^\top F\rvW\mid x,u\right)=\lambda^\top FD_\theta(\lambda)\left(\bar w+A_wx+B_wu\right),
\end{equation}
where $D_\theta(\lambda)$ is diagonal. Writing $z:=F^\top\lambda$, its diagonal entries are
\begin{equation}
\label{eq:polyhedral_risk_scaling}
D_\theta(\lambda)^{(ii)}:=
\begin{cases}
\displaystyle\frac{2}{1+\sqrt{1-\theta z^{(i)}/\kappa^{(i)}}}, & \text{IG},\\
\displaystyle\frac{e^{\theta z^{(i)}}-1}{\theta z^{(i)}}, & \text{Poisson},
\end{cases}
\end{equation}
with the ratio in the Poisson case defined as $1$ when $\theta z^{(i)}=0$. Note that $D_0(\lambda)=I$ for both models. For IG disturbances, $\mathcal D_\theta=\{\lambda\in\mathbb R_{\geq0}^{n_x}:\theta F^\top\lambda\leq\kappa\}$ when $\theta>0$, with $\kappa=[\kappa^{(1)}\cdots\kappa^{(n_d)}]^\top$, while $\mathcal D_\theta=\mathbb R_{\geq0}^{n_x}$ when $\theta\leq0$. For Poisson disturbances, $\mathcal D_\theta=\mathbb R_{\geq0}^{n_x}$ for every $\theta\in\mathbb R$.

Let $\mu_\theta(\lambda):=r+\bigl(B+FD_\theta(\lambda)B_w\bigr)^\top\lambda\in\mathbb R^{n_u}$, and let $E_i$ denote the $i$th row of $E$. Assumption~\ref{ass:affine_bellman} then holds for both models with\footnote{The minimization over $u$ follows that in \cite{Rantzer_2022}.}
\begin{subequations}
\label{eq:polyhedral_assumption_maps}
\begin{align}
&\Phi_\theta(\lambda)=s+\bigl(A+FD_\theta(\lambda)A_w\bigr)^\top\lambda-E^\top\bigl|\mu_\theta(\lambda)\bigr|,\\
&\Gamma_\theta(\lambda)=\lambda^\top FD_\theta(\lambda)\bar w,
\end{align}
\end{subequations}
and a measurable decision rule $a(\cdot;\lambda)$ attains the infimum in \eqref{eq:bellman_affinity} if and only if, for every $x\in\mathcal X$ and $i\in\{1,\ldots,n_u\}$, $a^{(i)}(x;\lambda)=-\operatorname{sgn}\bigl(\mu_\theta^{(i)}(\lambda)\bigr)E_ix$ when $\mu_\theta^{(i)}(\lambda)\neq0$, and $a^{(i)}(x;\lambda)\in[-E_ix,E_ix]$ otherwise.
\end{example}
The following example is concerned with systems evolving on the second-order cone.
\begin{example}
\label{ex:lorentz_cone}
Let $x=\left[x^{(1)},\tilde x^\top\right]^\top\in\mathbb{R}^{n_x}$ with $\tilde x\in\mathbb R^{n_x-1}$. Let $\mathcal X=\{x\in\mathbb R^{n_x}:x^{(1)}\geq\|\tilde x\|\}$ be the second-order cone, which is self-dual, i.e., $\mathcal X^\ast=\mathcal X$. Set $n_u=n_x$ and $\mathcal U(x)=\{u\in\mathcal X:u^{(1)}\leq\gamma x^{(1)}\}$, where $\gamma>0$, and consider
\begin{equation}
\label{eq:lorentz_model}
f(x,u,w)=
\begin{bmatrix}
\alpha x^{(1)}-u^{(1)}\\
A\tilde x-\tilde u
\end{bmatrix}+w,~
\ell(x,u,w)=s^\top x+r^\top u,
\end{equation}
where $A\in\mathbb R^{(n_x-1)\times(n_x-1)}$ satisfies $\|A\|+2\gamma\leq\alpha$, $s\in\mathcal{X}$ such that $s^{(1)}\geq\|\tilde s\|+\gamma\max\{\|\tilde r\|-r^{(1)},0\}$ (so that $s^\top x+r^\top u\geq0$ for all $x\in\mathcal{X}$ and $u\in\mathcal{U}(x)$).

Consider the disturbance
\begin{equation}
\label{eq:lorentz_disturbance}
\rvW=\sum_{k=1}^{\rvN}
\begin{bmatrix}
1\\
\rvV_k
\end{bmatrix},~
\rvN\mid(x,u)\sim\operatorname{Poisson}(\bar w+h^\top x),
\end{equation}
where $\bar w>0$ and $h\in\mathcal X$. Conditionally on $(x,u)$, the variables $\rvN,\rvV_1,\rvV_2,\ldots,\rvV_k$ are independent, and the $\rvV_1,\ldots,\rvV_k$ are identically distributed with $\mathbb P(\rvV_k=e_i)=\mathbb P(\rvV_k=-e_i)=1/(2(n_x-1))$ for $i\in\{1,\ldots,n_x-1\}$, where $e_i$ is the $i$th standard basis vector of $\mathbb R^{n_x-1}$. Since each arrival belongs to $\mathcal X$, we have $\rvW\in\mathcal X$. Moreover, $\|A\tilde x-\tilde u\|\leq\|A\|x^{(1)}+u^{(1)}\leq\alpha x^{(1)}-u^{(1)}$, and hence $\mathcal A(x)=\mathcal U(x)$ for every $x\in\mathcal X$.
For $\lambda=\left[\lambda^{(1)},\tilde{\lambda}^\top\right]^\top\in\mathcal X$, define
\begin{equation}
\label{eq:lorentz_arrival_risk}
\nu_{\theta,w}(\lambda):=
\begin{cases}
\displaystyle\frac{1}{\theta}\left(
\frac{e^{\theta\lambda^{(1)}}}{n_x-1}
\sum_{i=1}^{n_x-1}\cosh(\theta\tilde\lambda^{(i)})-1
\right),&\theta\neq0,\\[2ex]
\lambda^{(1)},&\theta=0.
\end{cases}
\end{equation}
Here, we have $\mathcal D_\theta=\mathcal X$ for every $\theta\in\mathbb R$, and, for any $\lambda\in\mathcal{X}^\ast$, $\mathcal R_\theta(\lambda^\top\rvW\mid x,u)=\nu_{\theta,w}(\lambda)(\bar w+h^\top x)$.

Assumption~\ref{ass:affine_bellman} then holds with
\begin{subequations}
\label{eq:lorentz_maps}
\begin{align}
&\begin{multlined}
\Phi_\theta(\lambda)=s+\nu_{\theta,w}(\lambda)h\\
+\begin{bmatrix}
\alpha\lambda^{(1)}-\gamma\max\{\lambda^{(1)}-r^{(1)}+\|\tilde\lambda-\tilde r\|,0\}\\
A^\top\tilde\lambda
\end{bmatrix},
\end{multlined}\\
&\Gamma_\theta(\lambda)=\bar w\,\nu_{\theta,w}(\lambda).
\end{align}
\end{subequations}
The minimum in \eqref{eq:opt_sol} is attained by $a(x;\lambda)=\left[a^{(1)}(x;\lambda) ~\tilde{a}^{\top}(x;\lambda)\right]^\top$ with $a^{(1)}(x;\lambda)=\gamma x^{(1)}$ when $\lambda^{(1)}+\|\tilde\lambda-\tilde r\|>r^{(1)}$, and $a^{(1)}(x;\lambda)=0$ otherwise. Moreover, $\tilde a(x;\lambda)=a^{(1)}(x;\lambda)(\tilde\lambda-\tilde r)/\|\tilde\lambda-\tilde r\|$ when $\tilde\lambda\neq\tilde r$, and $\tilde a(x;\lambda)=0$ otherwise. For $n_x=4$, the state can be interpreted as the Stokes vector of partially polarized light \cite{Goldstein_2011}, where $x^{(1)}$ is the total intensity and $\tilde x$ describes polarization.
\end{example}

\begin{remark}
\label{rem:classexample}
Examples~\ref{ex:polyhedral_ex1intro} and~\ref{ex:lorentz_cone} are instances of problems with the structure
\begin{equation}
\label{eq:condstruct1}
\mathcal{R}_\theta\left(\mathcal{H}(x,u,\rvW,\lambda)\mid x,u\right)
=\alpha_\theta(\lambda)+\eta_\theta^\top(\lambda)x+\mu_\theta^\top(\lambda)u,
\end{equation}
for all $\lambda\in\mathcal{D}_\theta$, $x\in\mathcal{X}$, and $u\in\mathcal{A}(x)$. In both examples, $\mathcal{P}:=\{(x,u):x\in\mathcal{X},\ u\in\mathcal{A}(x)\}$ is a proper cone, and \eqref{eq:bellman_affinity} holds if and only if there exists $\tilde{\eta}_\theta:\mathcal{D}_\theta\to\mathbb{R}^{n_x}$ such that
\begin{equation}
\label{eq:linear_cond1}
\inf_{u\in\mathcal{A}(x)} \mu_\theta^\top(\lambda)u=\tilde{\eta}_\theta^\top(\lambda)x,~ x\in\mathcal{X}.
\end{equation}
Conditions on $\mathcal{P}$ under which \eqref{eq:linear_cond1} holds are given in \cite[Theorem~2]{Pates_Rantzer_2024}.
\end{remark}

The following example does not satesfy \eqref{eq:condstruct1}.
\begin{example}
\label{ex:poisson_assumption}
Let $\mathcal X=\mathbb R_{\geq0}^{n_x}$ and $\mathcal U(x)=\mathbb R^{n_x}$, and consider
\begin{equation}
\label{eq:poisson_assumption_model}
f(x,u,w)=b+u+w,~\ell(x,u,w)=s^\top x+\frac{1}{2}u^\top Ru,
\end{equation}
where $b\in\mathbb R_{>0}^{n_x}$, $s\geq0$, and $R=\operatorname{diag}(r_1,\ldots,r_{n_x})$ with $r_1,\ldots,r_{n_x}>0$. Consider $\rvW\in\mathbb R_{\geq0}^{n_x}$ with conditionally independent components such that
\begin{equation}
\label{eq:poisson_assumption_distribution}
\rvW^{(i)}\mid(x,u)\sim\operatorname{Poisson}\left(\bar w^{(i)}+(Cx)^{(i)}\right),
\end{equation}
for $i\in\{1,\ldots,n_x\}$, where $\bar w\geq0$ and $C\geq0$. The mean of this term consists of an exogenous part $\bar w$ and a part $Cx$ generated by the current state, in which each unit of component $j$ generates on average $C^{(ij)}$ units of component $i$. The dynamics here are related to multitype branching processes with immigration \cite{Boxma_Ivanovs_Kosinski_Mandjes_2011}. These processes can be used to model computing workloads that generate additional tasks or population growth through reproduction. Since $\rvW=0$ occurs with positive probability, state viability gives $\mathcal A(x)=\{u\in\mathbb R^{n_x}:u\geq-b\}$. Moreover, $\mathcal D_\theta=\mathbb R_{\geq0}^{n_x}$ for every $\theta\in\mathbb R$. For $\lambda\in\mathcal D_\theta$, define $\nu_{\theta,w}(\lambda)\in\mathbb R_{\geq0}^{n_x}$ componentwise by $\nu_{\theta,w}^{(i)}(\lambda)=(e^{\theta\lambda^{(i)}}-1)/\theta$ for $\theta\neq0$, and $\nu_{0,w}(\lambda)=\lambda$. Then
\begin{equation}
\label{eq:poisson_ex2_risk}
\mathcal R_\theta\left(\lambda^\top\rvW\mid x,u\right)=\nu_{\theta,w}^\top(\lambda)(\bar w+Cx),
\end{equation}
and Assumption~\ref{ass:affine_bellman} holds since
\begin{subequations}
\label{eq:poisson_assumption_maps}
\begin{align}
&\Phi_\theta(\lambda)=s+C^\top \nu_{\theta,w}(\lambda),\\
&\Gamma_\theta(\lambda)=\bar w^\top \nu_{\theta,w}(\lambda)+\lambda^\top b+\min_{u\geq-b}\left\{\frac{1}{2}u^\top Ru+\lambda^\top u\right\}.
\end{align}
\end{subequations}
The minimum in \eqref{eq:opt_sol} is attained by $a^{(i)}(x;\lambda)=\max\{-b^{(i)},-\lambda^{(i)}/r_i\}$ for $i\in\{1,\ldots,n_x\}$.
\end{example}
\begin{remark}
\label{rem:risk_neutral_first_moment}
In Examples~\ref{ex:polyhedral_ex1intro}, \ref{ex:lorentz_cone}, and~\ref{ex:poisson_assumption}, $f$ and $\ell$ are affine in $w$. For $\theta=0$, linearity of expectation therefore gives $\mathcal R_0(\mathcal H(x,u,\rvW,\lambda)\mid x,u)=\mathcal H(x,u,\mathbb E[\rvW\!\mid \! x,u],\lambda)$. Hence, for a given admissible set $\mathcal A(x)$, the maps $\Phi_0$ and $\Gamma_0$ with the corresponding rule $a(\cdot;\lambda)$ depend on $\mathbb E[\rvW\mid x,u]$ only. For instance, any nonnegative disturbance with $\mathbb E[\rvW\mid x,u]=\bar w+A_wx+B_wu$ gives the same maps $\Phi_0$ and $\Gamma_0$ in Example~\ref{ex:polyhedral_ex1intro} ($D_0(\lambda)=I$).
\end{remark}
Note that the examples in \cite{Ohlin_Pates_Arcak_2025,semilineardynamicprogramminganalysis} also satesfy \ref{ass:affine_bellman}, including the postive bilinear system example (Example 2.1) in \cite{semilineardynamicprogramminganalysis}.

Additional examples, such as problems on PSD cones with LQG as a special case, will be discussed in Section~\ref{sec:examples}.

Finally, we state the following lemma which characterizes important properties of $\Gamma_\theta$ and $\Phi_\theta$ in \eqref{eq:bellman_affinity}.
\begin{lemma}
\label{lem:PhiGamma}
Let $\lambda\in\mathcal{D}_\theta$, and consider $\Gamma_\theta(\lambda)$ and $\Phi_\theta(\lambda)$ defined as in \eqref{eq:bellman_affinity}, then $\Gamma_\theta(\lambda)\geq0$ and $\Phi_\theta(\lambda)\in\mathcal X^\ast$. Furthermore, for $\lambda,\widetilde\lambda\in\mathcal D_\theta$ with $\lambda\geq_{\mathcal X^\ast}\widetilde\lambda$, we have that $\Gamma_\theta(\lambda)\geq\Gamma_\theta(\widetilde\lambda)$ and $\Phi_\theta(\lambda)\geq_{\mathcal X^\ast}\Phi_\theta(\widetilde\lambda)$, i.e., $\Gamma_\theta$ and $\Phi_\theta$ are order preserving with respect to the order on $\mathcal{X}^\ast$. If $\theta\leq0$, then $\Gamma_\theta$ is concave, and $\Phi_\theta$ is concave with respect to the order induced by $\mathcal X^\ast$, i.e., for $\lambda,\widetilde\lambda\in\mathcal X^\ast$ and $\tau\in[0,1]$,
\begin{equation}
\label{eq:Phi_order_concavity}
\Phi_\theta(\tau\lambda+(1-\tau)\widetilde\lambda)\geq_{\mathcal X^\ast}\tau\Phi_\theta(\lambda)+(1-\tau)\Phi_\theta(\widetilde\lambda).
\end{equation}
\end{lemma}

\begin{proof}
Combining \eqref{eq:Rtheta_geq0} and \eqref{eq:bellman_affinity} yields $\Gamma_\theta(\lambda)+\Phi_\theta(\lambda)^\top x\geq0$ for all $x\in\mathcal X$, which implies that $\Gamma_\theta(\lambda)\geq0$ and $\Phi_\theta(\lambda)\in\mathcal X^\ast$ by \eqref{eq:rep}.

To prove the order preserving properties, consider $\lambda,\widetilde\lambda\in\mathcal D_\theta$ with $\lambda\geq_{\mathcal X^\ast}\widetilde\lambda$. Then $\widetilde\lambda^\top x\leq\lambda^\top x$ for every $x\in\mathcal X$, so the monotonicity of $T_\theta$ in \eqref{eq:mono} and \eqref{eq:bellman_affinity} give $\Gamma_\theta(\widetilde\lambda)+\Phi_\theta(\widetilde\lambda)^\top x\leq\Gamma_\theta(\lambda)+\Phi_\theta(\lambda)^\top x$ for all $x\in\mathcal X$, which implies $\Gamma_\theta(\lambda)\geq\Gamma_\theta(\widetilde\lambda)$ and $\Phi_\theta(\lambda)\geq_{\mathcal X^\ast}\Phi_\theta(\widetilde\lambda)$ by \eqref{eq:rep}.

Now suppose $\theta\leq0$, so that $\mathcal D_\theta=\mathcal X^\ast$ by Proposition~\ref{prop:rs_domain}. Fix $\lambda,\widetilde\lambda\in\mathcal X^\ast$ and $\tau\in[0,1]$, and let $\lambda_\tau:=\tau\lambda+(1-\tau)\widetilde\lambda$. For every $x\in\mathcal X$ and $u\in\mathcal A(x)$, \eqref{eq:convexthetaR} divided by $\theta<0$, or \eqref{eq:expectationR} for $\theta=0$, shows that $\mathcal R_\theta(\mathcal H(x,u,\rvW,\lambda_\tau)\mid x,u)$ is at least $\tau\mathcal R_\theta(\mathcal H(x,u,\rvW,\lambda)\mid x,u)+(1-\tau)\mathcal R_\theta(\mathcal H(x,u,\rvW,\widetilde\lambda)\mid x,u)$. Taking the infimum over $u\in\mathcal A(x)$ and using \eqref{eq:bellman_affinity} gives
\begin{multline*}
\Gamma_\theta(\lambda_\tau)+\Phi_\theta(\lambda_\tau)^\top x\geq\\
\tau\Gamma_\theta(\lambda)+(1-\tau)\Gamma_\theta(\widetilde\lambda)+\left(\tau\Phi_\theta(\lambda)+(1-\tau)\Phi_\theta(\widetilde\lambda)\right)^\top x
\end{multline*}
for every $x\in\mathcal X$. By \eqref{eq:rep}, $\Gamma_\theta(\lambda_\tau)\geq\tau\Gamma_\theta(\lambda)+(1-\tau)\Gamma_\theta(\widetilde\lambda)$ and \eqref{eq:Phi_order_concavity} holds.
\end{proof}

\subsection{Finite-Horizon Problem}
\label{sec:finite_horizon}

Using Assumptions~\ref{ass:basic} and~\ref{ass:affine_bellman}, together with Lemmas~\ref{lem:PhiGamma} and~\ref{lem:finite_horizon_verification}, we obtain the following proposition.
\begin{proposition}
\label{prop:finite_horizon_affine}
Let $c_0=0$ and $\lambda_0\in\mathcal D_\theta$, and, for $N\geq0$, define
\begin{align}
\label{eq:affine_recursion}
\lambda_{N+1}&=\Phi_\theta(\lambda_N), &c_{N+1}&=c_N+\Gamma_\theta(\lambda_N),
\end{align}
whenever $\lambda_N\in\mathcal D_\theta$. Then, for every $N\geq0$ for which the recursion \eqref{eq:affine_recursion} is well defined, it holds that $c_N\geq0$, and the optimal value in \eqref{eq:rs_JN2} satisfies 
\begin{equation}
\label{eq:finite_horizon_affine_value}
J_{N,\lambda_0}(x)=c_N+\lambda_N^\top x,~x\in\mathcal X,
\end{equation}
{\color{black}
with $a_t^\star(x)$ in \eqref{eq:rs_optimal_feedback} given as 
\begin{equation}
\label{eq:finite_horizon_optimal_policy}
a_t^\star(x)=a(x;\lambda_{N-t-1}),~ t\in\{0,\ldots,N-1\},
\end{equation}
where $a(\cdot;\lambda)$ is defined in \eqref{eq:opt_sol}. 
}
Moreover, $c_{N+1}\geq c_N$ 
and if $\Phi_\theta(\lambda_0)\geq_{\mathcal X^\ast}\lambda_0$, then $\lambda_{N+1}\geq_{\mathcal X^\ast}\lambda_N$ and
\begin{equation}
\label{eq:finite_horizon_monotonicity}
J_{N+1,\lambda_0}(x)\geq J_{N,\lambda_0}(x),~ x\in\mathcal X,
\end{equation}
for every $N\geq0$. In particular, when $\lambda_0=0$, then \eqref{eq:finite_horizon_monotonicity} holds.
\end{proposition}
\begin{proof}
We consider only horizons for which the required iterates in
\eqref{eq:affine_recursion} are well defined. We first prove \eqref{eq:finite_horizon_affine_value} by induction on $N$. For $N=0$, the identity follows from
$J_{0,\lambda_0}(x)=\lambda_0^\top x$ and $c_0=0$. If it holds for some $N\geq0$ with $\lambda_N\in\mathcal D_\theta$, then \eqref{eq:rs_JN1},
translation invariance, and \eqref{eq:bellman_affinity} give
\begin{multline*}
J_{N+1,\lambda_0}(x)=(T_\theta J_{N,\lambda_0})(x)
\\=c_N+\Gamma_\theta(\lambda_N)
+\Phi_\theta(\lambda_N)^\top x
=c_{N+1}+\lambda_{N+1}^\top x.
\end{multline*}
This proves \eqref{eq:finite_horizon_affine_value} by induction.
Moreover, Lemma~\ref{lem:PhiGamma} gives $\lambda_N\in\mathcal{X}^\ast$ and $c_{N+1}-c_N=\Gamma_\theta(\lambda_N)\geq 0$, so that $c_N\geq0$ since $c_0=0$. Hence, $J_{k,\lambda_0}$ is finite, nonnegative, and Borel measurable for $k\in\{0,\ldots,N\}$, and Lemma~\ref{lem:finite_horizon_verification} gives \eqref{eq:rs_JN2} {\color{black} with \eqref{eq:rs_optimal_feedback} given by \eqref{eq:finite_horizon_optimal_policy} using translation invariance and Assumption~\ref{ass:affine_bellman}}.

Now suppose that $\Phi_\theta(\lambda_0)\geq_{\mathcal X^\ast}\lambda_0$. Then $\lambda_1\geq_{\mathcal X^\ast}\lambda_0$, and, by induction, the order preservation of $\Phi_\theta$ gives $\lambda_{N+1}=\Phi_\theta(\lambda_N)\geq_{\mathcal X^\ast}\Phi_\theta(\lambda_{N-1})=\lambda_N$ for all $N\geq1$. Consequently, $J_{N+1,\lambda_0}(x)-J_{N,\lambda_0}(x)=c_{N+1}-c_N+(\lambda_{N+1}-\lambda_N)^\top x\geq0$ for every $x\in\mathcal X$, which proves \eqref{eq:finite_horizon_monotonicity}. Finally, when $\lambda_0=0$, Lemma~\ref{lem:PhiGamma} gives $\Phi_\theta(0)\in\mathcal X^\ast$, and hence the same conclusion holds for $J_{N}$.
\end{proof}
\begin{remark}
\label{rem:conditions_welldefinedrec}
    By Lemma~\ref{lem:PhiGamma}, the recursion in \eqref{eq:affine_recursion} is well defined for every $N\geq0$ whenever $\mathcal D_\theta=\mathcal X^\ast$. In particular, this holds for $\theta\leq0$ by Proposition~\ref{prop:rs_domain}. When $\mathcal D_\theta\subsetneq\mathcal X^\ast$, as may occur for $\theta>0$ in Example~\ref{ex:polyhedral_ex1intro}, a sufficient condition is the existence of $\bar\lambda\in\mathcal D_\theta$ such that
\begin{equation}
\label{eq:recursion_supersolution}
\bar\lambda\geq_{\mathcal X^\ast}\Phi_\theta(\bar\lambda), 
\end{equation}
with $\bar\lambda\geq_{\mathcal X^\ast}\lambda_0$.
To see this, order preservation from Lemma~\ref{lem:PhiGamma} and downward closedness from Proposition~\ref{prop:rs_domain}, applied inductively to \eqref{eq:affine_recursion}, give $\lambda_N\in\mathcal D_\theta$ and $\bar\lambda\geq_{\mathcal X^\ast}\lambda_N$ for every $N\geq0$. The next subsection shows that the existence of $\bar\lambda\in\mathcal D_\theta$ satisfying \eqref{eq:recursion_supersolution} is equivalent to the existance of a fixed point of $\Phi_\theta$, under additional continuity assumptions on $\Phi_\theta$ and $\Gamma_\theta$.
\end{remark}

\subsection{Time-Average Infinite-Horizon Problem}
\label{sec:average_horizon}

In this subsection, we consider the time-average infinite-horizon problem in \eqref{eq:average_cost_ocp}. By Proposition~\ref{prop:finite_horizon_affine}, the relative values $J_N(x)-J_N(0)=\lambda_N^\top x$ obtained from \eqref{eq:affine_recursion} with $\lambda_0=0$ and $c_0=0$ are nondecreasing in $N$ whenever the recursion is well defined. Theorem~\ref{thm:average_cost} below show that they are bounded uniformly in $N$ if and only if $\Phi_\theta$ has a fixed point in $\mathcal D_\theta$. It also provides equivalent conditions for fixed-point existence and shows how the smallest fixed point characterizes the optimal time-average cost and a stationary optimal policy. One of these conditions yields an optimization characterization of the smallest fixed point (Corollary~\ref{cor:supersolution_program}).

Before presenting the infinite-horizon result, we make the following regularity assumption on $\Phi_\theta(\lambda)$ and $\Gamma_\theta(\lambda)$.

\begin{assumption}
\label{ass:rs_continuity}
The maps $\Phi_\theta:\mathcal D_\theta\to\mathbb R^{n_x}$ and $\Gamma_\theta:\mathcal D_\theta\to\mathbb R$ defined in \eqref{eq:bellman_affinity} are continuous.
\end{assumption}

We now present the main theorem for the infinite-horizon problem in \eqref{eq:average_cost_ocp}.
\begin{theorem}
\label{thm:average_cost}
Consider the recursion in \eqref{eq:affine_recursion} with $\lambda_0=0$ and $c_0=0$. Under Assumption~\ref{ass:rs_continuity}, the following statements are equivalent.
\begin{enumerate}
\item[(i)] For every $N\geq0$, the value function $J_N(x)$ is finite for every $x\in\mathcal X$, and there exists $\bar\lambda\in\mathcal D_\theta$ such that, for every $x\in\mathcal X$, $\sup_{N\geq0}\bigl(J_N(x)-J_N(0)\bigr)\leq\bar\lambda^\top x$.
\item[(ii)] The map $\Phi_\theta$ has a fixed point in $\mathcal D_\theta$.
\item[(iii)] There exists $\bar\lambda\in\mathcal D_\theta$ satisfying \eqref{eq:recursion_supersolution}.
\item[(iv)] There exist $\bar\lambda\in\mathcal D_\theta$, $\bar\rho\geq0$, and a stationary policy $\bar\pi=\{\bar a\}_{t\in\mathbb N}\in\Pi$ such that, for every $N\geq1$ and $x\in\mathcal X$,
\begin{equation}
\label{eq:finite_step_stabilizability}\mathcal R_\theta\left(S_{N,\bar\pi}^{\bar\lambda}(x)\right)\leq N\bar\rho+\bar\lambda^\top x.
\end{equation}
\end{enumerate}
If these equivalent conditions hold, then the recursion in \eqref{eq:affine_recursion} is well defined for every $N\geq0$ and $\lambda_N\to\lambda_\star$, where $\lambda_\star\in\mathcal D_\theta$ is a fixed point of $\Phi_\theta$ such that $\bar\lambda\geq_{\mathcal X^\ast}\lambda_\star$ for every $\bar\lambda\in\mathcal D_\theta$ satisfying \eqref{eq:recursion_supersolution}. In particular, $\lambda_\star$ is the smallest fixed point of $\Phi_\theta$ in $\mathcal D_\theta$ with respect to the order induced by $\mathcal X^\ast$. Moreover, for every $x\in\mathcal X$, $\sup_{N\geq0}\bigl(J_N(x)-J_N(0)\bigr)=\lambda_\star^\top x$.

With $\rho_\star:=\Gamma_\theta(\lambda_\star)$, the optimal average cost in \eqref{eq:average_cost_ocp} is
\begin{equation}
\label{eq:avg_rho}
J_{\mathrm{av}}^\star(x_0)=\rho_\star,~x_0\in\mathcal X.
\end{equation}
Furthermore, the stationary feedback policy $\pi_\star:=\{a(\cdot;\lambda_\star)\}_{t\in\mathbb N}$, given by
\begin{equation}
\label{eq:opt_sol_star}
\rvU_t=a(\rvX_t;\lambda_\star),~t\in\mathbb N,
\end{equation}
with $a(\cdot;\lambda_\star)$ from \eqref{eq:opt_sol}, is globally average-cost optimal.
\end{theorem}

\begin{proof}
Whenever the recursion in \eqref{eq:affine_recursion} with $\lambda_0=0$ is well defined up to horizon $N$, Proposition~\ref{prop:finite_horizon_affine} gives $J_N(x)=c_N+\lambda_N^\top x$, $\lambda_N\in\mathcal X^\ast$, and, for $N\geq1$, $\lambda_N\geq_{\mathcal X^\ast}\lambda_{N-1}$. Hence,
\begin{equation}
\label{eq:relative_value}
J_N(x)-J_N(0)=\lambda_N^\top x,~x\in\mathcal X.
\end{equation}
Moreover, $\mathcal D_\theta$ is downward closed with respect to the order induced by $\mathcal X^\ast$ by Proposition~\ref{prop:rs_domain}. Consequently, for every $\bar\lambda\in\mathcal D_\theta$,
\begin{multline}
\label{eq:domination_welldefined}
\bar\lambda\geq_{\mathcal X^\ast}\lambda_N\text{ whenever }\lambda_N\text{ is defined}\\
\Longrightarrow\quad\lambda_N\in\mathcal D_\theta\text{ for every }N\geq0.
\end{multline}
To see this, we have $\lambda_0=0$ is defined, and if $\lambda_N$ is defined, then $\lambda_N\in\mathcal X^\ast$ and downward closedness give $\lambda_N\in\mathcal D_\theta$, so that $\lambda_{N+1}=\Phi_\theta(\lambda_N)$ is defined. In particular, the recursion in \eqref{eq:affine_recursion} is then well defined for every $N\geq0$.

We first prove ${\rm (i)}\Rightarrow{\rm (ii)}$. Let $\bar\lambda$ be as in condition~${\rm (i)}$. By \eqref{eq:relative_value}, $\lambda_N^\top x\leq\bar\lambda^\top x$ for every $x\in\mathcal X$, i.e., $\bar\lambda\geq_{\mathcal X^\ast}\lambda_N$, whenever $\lambda_N$ is defined. Hence, by \eqref{eq:domination_welldefined}, $\lambda_N\in\mathcal D_\theta$ for every $N\geq0$, and, for every $x\in\mathcal X$, the sequence $\{\lambda_N^\top x\}_{N\geq0}$ is nondecreasing and bounded above by $\bar\lambda^\top x$, and therefore converges. By \eqref{eq:cone_span_bidual}, $\mathcal X$ contains a basis $x_1,\ldots,x_{n_x}$ of $\mathbb R^{n_x}$, and convergence of $\{\lambda_N^\top x_i\}_{N\geq0}$ for every $i\in\{1,\ldots,n_x\}$ gives $\lambda_N\to\lambda_\star$ for some $\lambda_\star\in\mathbb R^{n_x}$. Since $\mathcal X^\ast$ is closed, passing to the limit in $\lambda_N\in\mathcal X^\ast$ and $\bar\lambda\geq_{\mathcal X^\ast}\lambda_N$ gives $\lambda_\star\in\mathcal X^\ast$ and $\bar\lambda\geq_{\mathcal X^\ast}\lambda_\star$, and hence $\lambda_\star\in\mathcal D_\theta$ by downward closedness. Finally, Assumption~\ref{ass:rs_continuity} gives $
\lambda_\star=\lim_{N\to\infty}\lambda_{N+1}=\lim_{N\to\infty}\Phi_\theta(\lambda_N)=\Phi_\theta(\lambda_\star)$, which proves condition~${\rm (ii)}$.

The implication ${\rm (ii)}\Rightarrow{\rm (iii)}$ is immediate, since every fixed point $\lambda\in\mathcal D_\theta$ of $\Phi_\theta$ satisfies $\lambda\geq_{\mathcal X^\ast}\Phi_\theta(\lambda)$ with equality.

We next prove ${\rm (iii)}\Rightarrow{\rm (i)}$. Let $\bar\lambda\in\mathcal D_\theta$ satisfy \eqref{eq:recursion_supersolution}. Since $\bar\lambda\geq_{\mathcal X^\ast}0=\lambda_0$, Remark~\ref{rem:conditions_welldefinedrec} shows that the recursion in \eqref{eq:affine_recursion} is well defined and $\bar\lambda\geq_{\mathcal X^\ast}\lambda_N$ for every $N\geq0$. Hence, Proposition~\ref{prop:finite_horizon_affine} gives finite $J_N(x)$, and \eqref{eq:relative_value} gives $J_N(x)-J_N(0)=\lambda_N^\top x\leq\bar\lambda^\top x$, for every $x\in\mathcal X$ and $N\geq0$, which is condition~${\rm (i)}$.

We next prove ${\rm (iii)}\Rightarrow{\rm (iv)}$. Let $\bar\lambda\in\mathcal D_\theta$ satisfy \eqref{eq:recursion_supersolution}, and set $\bar a(x)=a(x;\bar\lambda)$, $\bar\pi=\{\bar a\}_{t\in\mathbb N}\in\Pi$, and $\bar\rho=\Gamma_\theta(\bar\lambda)\geq0$, using Assumption~\ref{ass:affine_bellman} and Lemma~\ref{lem:PhiGamma}. We now prove \eqref{eq:finite_step_stabilizability} by induction. For $N=1$, \eqref{eq:bellman_affinity} and \eqref{eq:recursion_supersolution} give
$\mathcal R_\theta\left(S_{1,\bar\pi}^{\bar\lambda}(x)\right)=\Gamma_\theta(\bar\lambda)+\Phi_\theta(\bar\lambda)^\top x\leq\bar\rho+\bar\lambda^\top x$.
Now assume that \eqref{eq:finite_step_stabilizability} holds for $N$. Then by \eqref{eq:disturbance}, \eqref{eq:recursive_R}, \eqref{eq:monoticity_of_R}, we have
\begin{multline*}
\mathcal R_\theta\left(S_{N+1,\bar\pi}^{\bar\lambda}(x)\right)\leq\mathcal R_\theta\left(S_{N,\bar\pi}^{0}(x)+\bar\rho+\bar\lambda^\top\rvX_N\right)=\\\bar\rho+\mathcal R_\theta\left(S_{N,\bar\pi}^{\bar\lambda}(x)\right)\leq(N+1)\bar\rho+\bar\lambda^\top x.
\end{multline*}
Hence, ${\rm (iv)}$ hold by induction.

We now prove ${\rm (iv)}\Rightarrow{\rm (iii)}$. For $N=1$ and using \eqref{eq:bellman_affinity}, we have
$\Gamma_\theta(\bar\lambda)+\Phi_\theta(\bar\lambda)^\top x\leq\mathcal R_\theta\left(S_{1,\bar\pi}^{\bar\lambda}(x)\right)\leq\bar\rho+\bar\lambda^\top x$.
By \eqref{eq:rep}, we have $\bar\lambda\geq_{\mathcal X^\ast}\Phi_\theta(\bar\lambda)$. Therefore, conditions~${\rm (i)}$-${\rm (iv)}$ are equivalent.

We next establish the convergence and minimality conclusions. Under the equivalent conditions, the proof of ${\rm (i)}\Rightarrow{\rm (ii)}$ shows that the recursion in \eqref{eq:affine_recursion} is well defined for every $N\geq0$ and that $\lambda_N\to\lambda_\star$ for some $\lambda_\star\in\mathcal D_\theta$ with $\lambda_\star=\Phi_\theta(\lambda_\star)$. If $\bar\lambda\in\mathcal D_\theta$ satisfies \eqref{eq:recursion_supersolution}, which holds, in particular, for every fixed point of $\Phi_\theta$ in $\mathcal D_\theta$, then the proof of ${\rm (iii)}\Rightarrow{\rm (i)}$ gives $\bar\lambda\geq_{\mathcal X^\ast}\lambda_N$ for every $N\geq0$, and passing to the limit gives $\bar\lambda\geq_{\mathcal X^\ast}\lambda_\star$. Hence $\lambda_\star$ is the smallest fixed point of $\Phi_\theta$ in $\mathcal D_\theta$, and, more generally, the smallest element of $\mathcal D_\theta$ satisfying \eqref{eq:recursion_supersolution}. Moreover, since $\lambda_N^\top x\nearrow\lambda_\star^\top x$ for every $x\in\mathcal X$, \eqref{eq:relative_value} gives $
\sup_{N\geq0}\bigl(J_N(x)-J_N(0)\bigr)=\lambda_\star^\top x$.

We now prove the average-cost conclusions. Set $\rho_\star:=\Gamma_\theta(\lambda_\star)$. By \eqref{eq:affine_recursion}, $c_N=\sum_{k=0}^{N-1}\Gamma_\theta(\lambda_k)$, and since $\Gamma_\theta(\lambda_k)\to\rho_\star$ by Assumption~\ref{ass:rs_continuity}, Cesàro's theorem gives $c_N/N\to\rho_\star$. Also, $\lambda_N^\top x_0/N\to0$ for every $x_0\in\mathcal X$. For every $\pi\in\Pi$, \eqref{eq:rs_JN2} and \eqref{eq:finite_horizon_affine_value} give
$ J_N^\pi(x_0)\geq J_N(x_0)=c_N+\lambda_N^\top x_0$.
Dividing by $N$ and taking the upper limit gives $ \limsup_{N\to\infty}\frac{1}{N}J_N^\pi(x_0)\geq\rho_\star$, and hence $J_{\mathrm{av}}^\star(x_0)\geq\rho_\star$ by \eqref{eq:average_cost_ocp}.

We next prove the reverse inequality. Let $J_\star(x)=\lambda_\star^\top x$. Since $\Phi_\theta(\lambda_\star)=\lambda_\star$ and $\rho_\star=\Gamma_\theta(\lambda_\star)$, Assumption~\ref{ass:affine_bellman} gives
\begin{equation}
\label{eq:fixed_point_bellman}
(T_\theta J_\star)(x)=\rho_\star+J_\star(x),~x\in\mathcal X,
\end{equation}
where the minimum is attained by $a(\cdot;\lambda_\star)$ from \eqref{eq:opt_sol}. Accordingly, every iterate of \eqref{eq:affine_recursion} started from $\lambda_\star$ equals $\lambda_\star$, with $c_N=N\rho_\star$, and Proposition~\ref{prop:finite_horizon_affine}, applied with $\lambda_0=\lambda_\star$, shows that the optimal policy in \eqref{eq:finite_horizon_optimal_policy} is $\pi_{\star,N}:=\{a(\cdot;\lambda_\star)\}_{t=0}^{N-1}\in\Pi_N$, which is the restriction of the stationary policy in \eqref{eq:opt_sol_star}, and that
\begin{equation}
\label{eq:terminal_identity}
\mathcal R_\theta\left(S^{\lambda_\star}_{N,\pi_{\star,N}}(x_0)\right)=J_{N,\lambda_\star}(x_0)=N\rho_\star+\lambda_\star^\top x_0.
\end{equation}
Since the policy in \eqref{eq:opt_sol_star} is admissible and $\lambda_\star\in\mathcal X^\ast$, we have $\lambda_\star^\top\rvX_N\geq0$. Hence, \eqref{eq:Nstage_cost_terminal}, monotonicity in \eqref{eq:monoticity_of_R}, and \eqref{eq:terminal_identity} imply
\begin{equation}
\label{eq:stationary_upper_bound}
J_N^{\pi_\star}(x_0)\leq N\rho_\star+\lambda_\star^\top x_0.
\end{equation}
Dividing \eqref{eq:stationary_upper_bound} by $N$ and taking the upper limit gives $
\limsup_{N\to\infty}\frac{1 }{N}J_N^{\pi_\star}(x_0)\leq\rho_\star\leq J_{\mathrm{av}}^\star(x_0)$, which, together with \eqref{eq:average_cost_ocp}, proves \eqref{eq:avg_rho} and shows that the stationary policy in \eqref{eq:opt_sol_star} is globally average-cost optimal.
\end{proof}
Under Assumption~\ref{ass:rs_continuity} and the equivalent conditions of Theorem~\ref{thm:average_cost}, $J_N(x)-J_N(0)=\lambda_N^\top x\to\lambda_\star^\top x$ and $J_N(x)/N\to\rho_\star$ for every $x\in\mathcal X$ as proven above. Thus, $\lambda_\star^\top x$ describes the limiting dependence of the finite-horizon value on the initial state relative to $0$, while $\rho_\star$ is the optimal cost per stage. Moreover, \eqref{eq:terminal_identity} shows that the stationary policy in \eqref{eq:opt_sol_star} is optimal for every finite horizon when the terminal cost is $\lambda_\star^\top x$. 
Note that the theorem does not require uniqueness of the fixed point, and the recursion initialized at zero selects the smallest one.

\begin{remark}
\label{rem:one_step_horizons}
Taking $N=1$ in \eqref{eq:finite_step_stabilizability} gives
\begin{equation}
\label{eq:stabilizability}
\mathcal R_\theta\left(\mathcal H(x,\bar a(x),\rvW,\bar\lambda)\mid x,\bar a(x)\right)\leq\bar\rho+\bar\lambda^\top x,~x\in\mathcal X.
\end{equation}
For $\theta=0$, this is an expected dissipation inequality for the linear function $x\mapsto\bar\lambda^\top x$, with a risk-sensitive counterpart for $\theta\neq0$. Moreover, \eqref{eq:finite_step_stabilizability} provides a finite-horizon bound on the risk-sensitive evaluation of the accumulated stage cost and terminal cost under stationary policies satisfying condition~${\rm (iv)}$, including the optimal policy in \eqref{eq:opt_sol_star}. In Section~\ref{sec:detectability_performance}, we combine \eqref{eq:finite_step_stabilizability} with an additional observability assumption to establish uniform state bounds for $\theta\geq0$.
\end{remark}

The smallest fixed point also admits the following optimization characterization.
\begin{corollary}
\label{cor:supersolution_program}
Suppose that Assumption~\ref{ass:rs_continuity} holds, let $\xi\in\operatorname{int}\mathcal X$, and consider
\begin{equation}
\label{eq:supersolution_program}\inf_{\lambda\in\mathcal D_\theta}\left\{\xi^\top\lambda:~\lambda\geq_{\mathcal X^\ast}\Phi_\theta(\lambda)\right\}.
\end{equation}
This problem is feasible if and only if the equivalent conditions in Theorem~\ref{thm:average_cost} hold, in which case $\lambda_\star$ is its unique minimizer. Moreovre, if $\mathcal{D}_\theta$ is convex and $\Phi_\theta$ is convex with respect to the order induced by $\mathcal{X}^\ast$, then \eqref{eq:supersolution_program} is convex.
\end{corollary}
\begin{proof}
Feasibility is precisely condition~${\rm (iii)}$ of Theorem~\ref{thm:average_cost}. When \eqref{eq:supersolution_program} is feasible, the Theorem ~\ref{thm:average_cost} gives feasibility of $\lambda_\star$ and $\bar\lambda\geq_{\mathcal{X}^\ast}\lambda_\star$ for every feasible $\bar\lambda$. Since $\xi\in\operatorname{int}\mathcal X$, \eqref{eq:cone_interior_norm} gives $\xi^\top(\bar\lambda-\lambda_\star)>0$ whenever $\bar\lambda\neq\lambda_\star$, which proves uniqueness.
Finally, under the stated convexity assumptions, for feasible $\lambda,\widetilde\lambda$ and $\tau\in[0,1]$ we have
\begin{multline*}
\Phi_\theta(\tau\lambda+(1-\tau)\widetilde\lambda)\leq_{\mathcal X^\ast}\tau\Phi_\theta(\lambda)+(1-\tau)\Phi_\theta(\widetilde\lambda)\\\leq_{\mathcal X^\ast}\tau\lambda+(1-\tau)\widetilde\lambda.
\end{multline*}
Thus, \eqref{eq:supersolution_program} is a convex optimization problem. 
\end{proof}


\begin{remark}
\label{rem:supersolution_program_examples}
In Example~\ref{ex:poisson_assumption}, $\mathcal D_\theta=\mathbb R_{\geq0}^{n_x}$ and $\Phi_\theta(\lambda)=s+C^\top \nu_{\theta,w}(\lambda)$ with $C\geq0$. Thus, \eqref{eq:supersolution_program} is convex for $\theta\geq0$. It is also a linear program when $\theta=0$ and can be written using exponential cone constraints when $\theta>0$. In general, for $\theta\leq0$, Lemma~\ref{lem:PhiGamma} gives order concavity of $\Phi_\theta$, so \eqref{eq:supersolution_program} need not be convex. In Section~\ref{sec:detectability_performance}, a different convex characterization for $\theta\leq0$ will be given under an additional observability condition.
\end{remark}

\begin{remark}
\label{rem:stationary_finite_horizon}
As no practical application can run indefinitely, control policies are typically applied over finite and often relatively long horizons. For some applications, the stationary policy in \eqref{eq:opt_sol_star} can be computed once and implemented without computing and storing the time-varying decision rules in \eqref{eq:finite_horizon_optimal_policy}. In such cases, define the \emph{finite-horizon suboptimality gap} of the stationary policy by $\Delta_N(x):=J_N^{\pi_\star}(x)-J_N(x)$ for $N\geq1$. By \eqref{eq:stationary_upper_bound} and Proposition~\ref{prop:finite_horizon_affine}, we have
\begin{equation*}
0\leq\frac{\Delta_N(x)}{N}
\leq\rho_\star+\frac{\lambda_\star^\top x}{N}-\frac{J_N(x)}{N}.
\end{equation*}
Since $J_N(x)/N\to\rho_\star$ by Theorem~\ref{thm:average_cost}, it follows that
\begin{equation}
\label{eq:stationary_per_stage_difference}
\lim_{N\to\infty}\frac{\Delta_N(x)}{N}=0,~x\in\mathcal X.
\end{equation}
Therefore, on average, the finite-horizon suboptimality gap converges to zero as $N\to\infty$.
\end{remark}

\subsection{Zero Average Cost and Total Cost}
\label{sec:stochastic_total_cost}
We now provide a proposition with a condition for when the infintie horizon problem $J_{\star}(x_0)=\inf_{\pi\in\Pi}\mathcal R_\theta\left(\sum_{t=0}^{\infty}\ell(\rvX_t,\rvU_t,\rvW_t)\right)$ is finite and has $\eqref{eq:opt_sol_star}$ as the optimal solution for any $x_0\in\mathcal{X}$.  
\begin{proposition}
\label{prop:stochastic_total_cost}
Suppose that Assumption~\ref{ass:rs_continuity} holds and there exists $\bar\lambda\in\mathcal D_\theta$ such that
\begin{equation}
\label{eq:total_cost_supersolution}
\bar\lambda\geq_{\mathcal X^\ast}\Phi_\theta(\bar\lambda),~ \Gamma_\theta(\bar\lambda)=0.
\end{equation}
Then, for every $x\in\mathcal X$,
\begin{equation}
\label{eq:stochastic_total_cost}
\inf_{\pi\in\Pi}\mathcal R_\theta\left(\sum_{t=0}^{\infty}\ell(\rvX_t,\rvU_t,\rvW_t)\right)
=\inf_{\pi\in\Pi}\sup_{N\geq1}J_N^\pi(x)=\lambda_\star^\top x,
\end{equation}
where $\lambda_\star$ is the smallest fixed point from Theorem~\ref{thm:average_cost}, and the infimum is attained by the stationary policy in \eqref{eq:opt_sol_star}.
\end{proposition}
\begin{proof}
By Theorem~\ref{thm:average_cost} and Lemma~\ref{lem:PhiGamma}, \eqref{eq:total_cost_supersolution} gives $0\leq\Gamma_\theta(\lambda_N)\leq\Gamma_\theta(\lambda_\star)\leq\Gamma_\theta(\bar\lambda)=0$. Hence, $\rho_\star=0$ and $c_N=0$ for every $N\geq0$. For every $\pi\in\Pi$ and $N\geq1$, Proposition~\ref{prop:finite_horizon_affine} gives $J_N^\pi(x)\geq J_N(x)=\lambda_N^\top x$, while \eqref{eq:stationary_upper_bound} gives $J_N^{\pi_\star}(x)\leq\lambda_\star^\top x$. Taking suprema over $N$ and using $\lambda_N\to\lambda_\star$, we obtain
$\inf_{\pi\in\Pi}\sup_{N\geq1}J_N^\pi(x)=\sup_{N\geq1}J_N^{\pi_\star}(x)=\lambda_\star^\top x$. Finally, nonnegativity of the stage costs gives $S_{N,\pi}^0(x)\nearrow\sum_{t=0}^{\infty}\ell(\rvX_t,\rvU_t,\rvW_t)$ for every $\pi\in\Pi$. Monotone convergence for $\theta\geq0$ and dominated convergence for $\theta<0$ therefore give $\mathcal R_\theta(\sum_{t=0}^{\infty}\ell(\rvX_t,\rvU_t,\rvW_t))=\sup_{N\geq1}J_N^\pi(x)$, which proves the claim.
\end{proof}

\section{Finite-Step Observability and Its Consequences}
\label{sec:detectability_performance}
We now introduce a finite-step cost observability condition for risk-sensitive control on proper cones, analogous to the finite-step positive-cost condition in \cite[Assumption~2.1(d)]{semilineardynamicprogramminganalysis}, to strengthen the conclusions of Theorem~\ref{thm:average_cost}. Specifically, for $\theta\leq0$, it gives uniqueness of the fixed point $\lambda_\star$, geometric convergence from arbitrary $\lambda_0\in\mathcal X^\ast$ (instead of $\lambda_0=0$ as in Theorem~\ref{thm:average_cost}), and a bound, uniform in the horizon, on the finite-horizon suboptimality gap of $\pi_\star$. For $\theta\geq0$, it provides uniform state bounds as discussed in Remark~\ref{rem:one_step_horizons}.
\subsection{Finite-Step Cost Observability}
\label{sec:cost_observability}

\begin{assumption}[Finite-step cost observability]
\label{ass:finite_step_detectability}
There exist $N_{\mathrm d}\in\mathbb N\setminus\{0\}$ and $\eta\in\operatorname{int}\mathcal X^\ast$ such that, for every $x\in\mathcal X$ and every $\pi_{N_{\mathrm d}}\in\Pi_{N_{\mathrm d}}$,
\begin{equation}
\label{eq:finite_step_detectability}
\mathcal R_{-|\theta|}\left(S_{N_{\mathrm d},\pi_{N_{\mathrm d}}}^{0}(x)\right)\geq\eta^\top x,
\end{equation}
where the entropic-risk evaluation uses the law induced by $x$, $\pi_{N_{\mathrm d}}$, and $\rcd$ in \eqref{eq:general_dynamics}. 
\end{assumption}

Since $\eta\in\operatorname{int}\mathcal X^\ast$, \eqref{eq:cone_interior_norm} implies that there exist $\delta>0$ such that $\eta^\top x\geq\delta\|x\|$ for every $x\in\mathcal X$. Thus, \eqref{eq:finite_step_detectability} gives
$\mathcal R_{-|\theta|}\left(S_{N_{\mathrm d},\pi_{N_{\mathrm d}}}^{0}(x)\right)\geq\eta^\top x\geq\delta\|x\|$.
We refer to this condition as cost observability because the accumulated stage cost detects the size of the initial state. In the deterministic case, it gives $S_{N_{\mathrm d},\pi_{N_{\mathrm d}}}^{0}(x)\geq\delta\|x\|$, so no admissible policy can keep the accumulated cost bounded as the initial state grows. For $\theta=0$, the lower bound holds in expectation.
For $\theta\neq0$, \eqref{eq:finite_step_detectability} is equivalent to
\begin{equation}
\label{eq:finite_step_observability_exponential}
\mathbb E\left[e^{-|\theta| S_{N_{\mathrm d},\pi_{N_{\mathrm d}}}^{0}(x)}\right]\leq e^{-|\theta|\eta^\top x}.
\end{equation}
Applying Markov's inequality gives, for every $q\in(0,1)$,
\begin{multline}
\label{eq:finite_step_observability_probability}
\mathbb P\left(S_{N_{\mathrm d},\pi_{N_{\mathrm d}}}^{0}(x)\leq q\eta^\top x\right)\\
\leq e^{-|\theta|(1-q)\eta^\top x}\leq e^{-|\theta|(1-q)\delta\|x\|}.
\end{multline}
Thus, small accumulated costs remain possible, but the probability of incurring less than a fixed fraction of $\eta^\top x$ has an upper bound that decreases exponentially with $\|x\|$, uniformly over admissible policies.

The following lemma provides a test for observability based on the reccursion in \eqref{eq:affine_recursion} with $\lambda_0=0$ and $c_0=0$, which is necessary and sufficient when $\theta\leq0$.
Related finite-horizon positivity results have been established for deterministic positive systems \cite{DynamicProgramming_YuchaoRantzer} and risk-neutral stochastic semilinear dynamic programming for positive systems \cite{semilineardynamicprogramminganalysis}.

\begin{lemma}
\label{lem:detectability_coefficients}
Let $N_{\mathrm d}\geq1$, and consider the recursion in \eqref{eq:affine_recursion} with $\lambda_0=0$ and $c_0=0$, assumed to be well defined up to horizon $N_{\mathrm d}$ if $\theta>0$. If Assumption~\ref{ass:finite_step_detectability} holds with $N_{\mathrm d}$, then
\begin{equation}
\label{eq:detectability_coefficients}
\lambda_{N_{\mathrm d}}\in\operatorname{int}\mathcal X^\ast,
\end{equation}
and $\lambda_N\in\operatorname{int}\mathcal X^\ast$ for every $N\geq N_{\mathrm d}$ for which $\lambda_N$ is defined. In particular, $\lambda_\star\in\operatorname{int}\mathcal X^\ast$ whenever the equivalent conditions in Theorem~\ref{thm:average_cost} hold. Conversely, if $\theta\leq0$ and \eqref{eq:detectability_coefficients} holds, then Assumption~\ref{ass:finite_step_detectability} holds with $\eta=\lambda_{N_{\mathrm d}}$.
\end{lemma}
\begin{proof}
For every nonnegative random variable $\rvY$, $\mathcal R_\theta(\rvY)\geq\mathcal R_{-|\theta|}(\rvY)$, with equality for $\theta\leq0$ and the inequality following from Jensen's inequality for $\theta>0$. Hence, taking the infimum over policies in \eqref{eq:finite_step_detectability} and using Proposition~\ref{prop:finite_horizon_affine} gives $c_{N_{\mathrm d}}+(\lambda_{N_{\mathrm d}}-\eta)^\top x\geq0$, for all $x\in\mathcal X$.
By \eqref{eq:rep}, $\lambda_{N_{\mathrm d}}\geq_{\mathcal X^\ast}\eta$, so \eqref{eq:cone_interior_addition} gives $\lambda_{N_{\mathrm d}}\in\operatorname{int}\mathcal X^\ast$. The conclusions for $\lambda_N$ and $\lambda_\star$ follow from \eqref{eq:cone_interior_addition}, the monotonicity of $\lambda_N$ in Proposition~\ref{prop:finite_horizon_affine}, its convergence to $\lambda_\star$ in Theorem~\ref{thm:average_cost}, and closedness of $\mathcal X^\ast$. Conversely, for $\theta\leq0$ and every $\pi_{N_{\mathrm d}}\in\Pi_{N_{\mathrm d}}$, \eqref{eq:rs_JN2} gives
$\mathcal R_{-|\theta|}\left(S_{N_{\mathrm d},\pi_{N_{\mathrm d}}}^{0}(x)\right)\geq J_{N_{\mathrm d}}(x)
=c_{N_{\mathrm d}}+\lambda_{N_{\mathrm d}}^\top x\geq\lambda_{N_{\mathrm d}}^\top x$.
This proves the converse with $\eta=\lambda_{N_{\mathrm d}}$.
\end{proof}

\begin{remark}
\label{rem:negative_risk_observability_test}
For $\theta>0$, suppose that Assumption~\ref{ass:affine_bellman} also hold at $-\theta$, and consider $\widetilde\lambda_0=0,~\widetilde\lambda_{N+1}=\Phi_{-\theta}(\widetilde\lambda_N)$. Then Lemma~\ref{lem:detectability_coefficients}, applied at $-\theta$, shows that Assumption~\ref{ass:finite_step_detectability} holds with horizon $N_{\mathrm d}$ if and only if $\widetilde\lambda_{N_{\mathrm d}}\in\operatorname{int}\mathcal X^\ast$. One can then take $\eta=\widetilde\lambda_{N_{\mathrm d}}$. 
\end{remark}

\begin{remark}
\label{rem:intpropercone_detect}
\label{rem:detectability_example1}
Suppose that $\mathcal P=\{(x,u):x\in\mathcal X,\ u\in\mathcal U(x)\}$ is a proper cone and that $\mathcal R_{-|\theta|}(\ell(x,u,\rvW)\mid x,u)\geq s^\top x+r^\top u$ for every $(x,u)\in\mathcal P$, with $(s,r)\in\operatorname{int}\mathcal P^\ast$. For any $\zeta\in\operatorname{int}\mathcal X^\ast$, \eqref{eq:cone_interior_comparison} gives $\delta>0$ such that $(s-\delta\zeta,r)\in\mathcal P^\ast$. Hence, $\mathcal R_{-|\theta|}(\ell(x,u,\rvW)\mid x,u)\geq\delta\zeta^\top x$, and Assumption~\ref{ass:finite_step_detectability} holds with $N_{\mathrm d}=1$ and $\eta=\delta\zeta\in\operatorname{int}\mathcal X^\ast$ by \eqref{eq:cone_interior_addition}. This is the counterpart of the deterministic condition in \cite{Pates_Rantzer_2024} for $\ell(x,u)=s^\top x+r^\top u$. In Examples~\ref{ex:polyhedral_ex1intro} and~\ref{ex:lorentz_cone}, where $\ell(x,u)=s^\top x+r^\top u$, the condition $(s,r)\in\operatorname{int}\mathcal P^\ast$ is therefore sufficient for every $\theta\in\mathbb R$. Specifically, for Example~\ref{ex:polyhedral_ex1intro}, this amounts to $s-E^\top|r|>0$, with $\eta=s-E^\top|r|$. In Example~\ref{ex:poisson_assumption}, $\ell(x,u)\geq s^\top x$, so Assumption~\ref{ass:finite_step_detectability} holds with $N_{\mathrm d}=1$ and $\eta=s$ whenever $s>0$.
\end{remark}

\subsection{Convergence and Finite-Horizon Performance}
\label{sec:convergence_performance}

For $\theta\leq0$, the following theorem establishes uniqueness of the fixed point $\lambda_\star$ in Theorem~\ref{thm:average_cost} (in contrast to just being the smallest one), quantifies convergence to it (from any initial $\lambda_0\in\mathcal{X}^\ast$ instead of only $\lambda_0=0$), and bounds the finite-horizon suboptimality gap of $\pi_\star$ uniformly in the horizon.
Under Assumption~\ref{ass:rs_continuity}, the equivalent conditions of Theorem~\ref{thm:average_cost}, and \eqref{eq:detectability_coefficients}, consider the recursion in \eqref{eq:affine_recursion} with $\lambda_0=0$ and $c_0=0$, and define
\begin{equation}
\label{eq:detectability_margin}
\beta:=\max\left\{b\in[0,1]:~\lambda_{N_{\mathrm d}}\geq_{\mathcal X^\ast}b\lambda_\star\right\}.
\end{equation}
By \eqref{eq:relative_value}, $\beta$ is the largest fraction of $\lambda_\star^\top x$ attained within $N_{\mathrm d}$ stages, uniformly over $x\in\mathcal X$. Below, $\Phi_\theta^N$ denotes the $N$-fold composition of $\Phi_\theta$, with $\Phi_\theta^0$ the identity, and we use the convention $0^0=1$.

\begin{theorem}
\label{thm:detectability_consequences}
Let $\theta\leq0$, and suppose that Assumption~\ref{ass:rs_continuity}, the equivalent conditions in Theorem~\ref{thm:average_cost}, and Assumption~\ref{ass:finite_step_detectability} hold. Let $\beta$ be given by \eqref{eq:detectability_margin}. Then $\beta\in(0,1]$ and the following statements hold.
\begin{enumerate}
\item[(i)] For every $\lambda_0\in\mathcal X^\ast$, let $\mu\geq1$ be such that $\mu\lambda_\star\geq_{\mathcal X^\ast}\lambda_0$. Then, for every $N\geq0$ and $x\in\mathcal X$,
\begin{multline}
\label{eq:global_convergence}
-(1-\beta)^{\lfloor N/N_{\mathrm d}\rfloor}\lambda_\star^\top x\leq\left(\Phi_\theta^N(\lambda_0)-\lambda_\star\right)^\top x\\
\leq(\mu-1)(1-\beta)^{\lfloor N/N_{\mathrm d}\rfloor}\lambda_\star^\top x.
\end{multline}
In particular, $\lambda_\star$ is the unique fixed point of $\Phi_\theta$ in $\mathcal X^\ast$, and $J_{N,\lambda_0}(x)-J_{N,\lambda_0}(0)$ converges to $\lambda_\star^\top x$ for every $x\in\mathcal X$.
\item[(ii)] For every $x\in\mathcal X$ and $N\geq1$, the finite-horizon suboptimality gap satisfies
\begin{equation}
\label{eq:bounded_regret}0\leq\Delta_N(x)\leq\frac{N_{\mathrm d}\rho_\star}{\beta}+(1-\beta)^{\lfloor N/N_{\mathrm d}\rfloor}\lambda_\star^\top x.
\end{equation}
\item[(iii)] For every $\xi\in\operatorname{int}\mathcal X$, the problem
\begin{equation}
\label{eq:subsolution_program}
\sup_{\lambda\in\mathcal X^\ast}\left\{\xi^\top\lambda:~\Phi_\theta(\lambda)\geq_{\mathcal X^\ast}\lambda\right\}
\end{equation}
is a convex optimization problem whose unique maximizer is $\lambda_\star$.
\end{enumerate}
\end{theorem}
\begin{proof}
For $\theta\leq0$, Proposition~\ref{prop:rs_domain} gives $\mathcal D_\theta=\mathcal X^\ast$, and Lemma~\ref{lem:PhiGamma} gives order preservation and order concavity of $\Phi_\theta$. Both properties are preserved under composition, so they hold for $\Phi_\theta^N$ for every $N\geq0$. By Lemma~\ref{lem:detectability_coefficients}, $\lambda_{N_{\mathrm d}}$ and $\lambda_\star$ belong to $\operatorname{int}\mathcal X^\ast$. Thus, \eqref{eq:cone_interior_comparison} gives $\lambda_{N_{\mathrm d}}\geq_{\mathcal X^\ast}b\lambda_\star$ for sufficiently small $b>0$, so $\beta>0$, and the maximum in \eqref{eq:detectability_margin} is attained by compactness. Moreover, \eqref{eq:cone_interior_comparison}, applied with $p=\lambda_\star$ and $v=\lambda_0$, gives $\mu\geq1$ such that $\mu\lambda_\star\geq_{\mathcal X^\ast}\lambda_0$.

We first prove (i). Set $t_k:=1-(1-\beta)^k$. We show by induction that $\lambda_{kN_{\mathrm d}}\geq_{\mathcal X^\ast}t_k\lambda_\star$. The claim holds for $k=0$. If it holds for $k$, order preservation, order concavity, and \eqref{eq:detectability_margin} give \begin{multline*}
\lambda_{(k+1)N_{\mathrm d}}\geq_{\mathcal X^\ast}\Phi_\theta^{N_{\mathrm d}}(t_k\lambda_\star)\\
\geq_{\mathcal X^\ast}t_k\lambda_\star+(1-t_k)\lambda_{N_{\mathrm d}}\geq_{\mathcal X^\ast}t_{k+1}\lambda_\star.
\end{multline*}
For any $N\geq0$, let $k:=\lfloor N/N_{\mathrm d}\rfloor$. Since $kN_{\mathrm d}\leq N$, the monotonicity of $\lambda_N$ from Proposition~\ref{prop:finite_horizon_affine} gives $\lambda_N\geq_{\mathcal X^\ast}\lambda_{kN_{\mathrm d}}\geq_{\mathcal X^\ast}t_k\lambda_\star$. Hence,
\begin{equation}
\label{eq:sandwich_below}
\lambda_N\geq_{\mathcal X^\ast}t_{\lfloor N/N_{\mathrm d}\rfloor}\lambda_\star\geq_{\mathcal X^\ast}0,~N\geq0.
\end{equation}
For $\lambda_0$ and $\mu$ such that $\mu\lambda_\star\geq_{\mathcal X^\ast}\lambda_0$, order preservation and order concavity also give
$\lambda_N\leq_{\mathcal X^\ast}\Phi_\theta^N(\lambda_0)\leq_{\mathcal X^\ast}\Phi_\theta^N(\mu\lambda_\star)
\leq_{\mathcal X^\ast}\mu\lambda_\star-(\mu-1)\lambda_N$,
where the last inequality follows by applying order concavity to $\lambda_\star=\mu^{-1}(\mu\lambda_\star)+(1-\mu^{-1})0$. Combining this with \eqref{eq:sandwich_below}, and setting $k:=\lfloor N/N_{\mathrm d}\rfloor$, gives
\begin{multline*}
(t_{\lfloor N/N_{\mathrm d}\rfloor}-1)\lambda_\star\leq_{\mathcal X^\ast}\Phi_\theta^N(\lambda_0)-\lambda_\star\\
\leq_{\mathcal X^\ast}(\mu-1)(1-t_{\lfloor N/N_{\mathrm d}\rfloor})\lambda_\star.
\end{multline*}
Using $1-t_k=(1-\beta)^k$, we get \eqref{eq:global_convergence}.
Since $\beta\in(0,1]$, \eqref{eq:global_convergence} gives $(\Phi_\theta^N(\lambda_0)-\lambda_\star)^\top x\to0$ for every $x\in\mathcal X$. By \eqref{eq:cone_span_bidual}, $\mathcal X$ contains a basis of $\mathbb R^{n_x}$, and therefore $\Phi_\theta^N(\lambda_0)\to\lambda_\star$. Now let $\widetilde\lambda\in\mathcal X^\ast$ be any fixed point of $\Phi_\theta$. Since $\Phi_\theta^N(\widetilde\lambda)=\widetilde\lambda$ for every $N\geq0$, the preceding convergence with $\lambda_0=\widetilde\lambda$ gives $\widetilde\lambda=\lambda_\star$, proving uniqueness. Finally, Proposition~\ref{prop:finite_horizon_affine} gives $J_{N,\lambda_0}(x)-J_{N,\lambda_0}(0)=\Phi_\theta^N(\lambda_0)^\top x\to\lambda_\star^\top x$ for every $x\in\mathcal X$.

We next prove (ii). Recall that $\Delta_N(x)=J^{\pi_\star}_{N}(x)-J_N(x)$ as defined in Remark~\ref{rem:stationary_finite_horizon}. By \eqref{eq:stationary_upper_bound} and Proposition~\ref{prop:finite_horizon_affine}, we have
\begin{equation}
\label{eq:diff_finitehorizon_stationary_nobounds}
0\leq J_N^{\pi_\star}(x)-J_N(x)
\leq N\rho_\star-c_N+(\lambda_\star-\lambda_N)^\top x.
\end{equation}
By \eqref{eq:sandwich_below}, we obtain
\begin{equation}
\label{eq:relative_bound_finitehorizon_stationary}
(\lambda_\star-\lambda_N)^\top x
\leq(1-\beta)^{\lfloor N/N_{\mathrm d}\rfloor}\lambda_\star^\top x.
\end{equation}
It remains to bound $N\rho_\star-c_N$. By Lemma~\ref{lem:PhiGamma}, $\Gamma_\theta$ is order preserving and concave, with $\Gamma_\theta(0)\geq0$. For every $j\geq0$, let $k:=\lfloor j/N_{\mathrm d}\rfloor$. By \eqref{eq:sandwich_below} and Theorem~\ref{thm:average_cost}, we have $t_k\lambda_\star\leq_{\mathcal X^\ast}\lambda_j\leq_{\mathcal X^\ast}\lambda_\star$. Hence,
$\rho_\star\geq\Gamma_\theta(\lambda_j)
\geq\Gamma_\theta(t_k\lambda_\star)
\geq t_k\rho_\star+(1-t_k)\Gamma_\theta(0)\geq t_k\rho_\star$.
Since $1-t_k=(1-\beta)^k$, it follows that $
0\leq\rho_\star-\Gamma_\theta(\lambda_j)\leq(1-\beta)^{\lfloor j/N_{\mathrm d}\rfloor}\rho_\star$.
Summing this inequality and using \eqref{eq:affine_recursion} gives
\begin{multline}
\label{eq:bias_bound_finitehorizon_stationary}
0\leq N\rho_\star-c_N
=\sum_{j=0}^{N-1}\left(\rho_\star-\Gamma_\theta(\lambda_j)\right)\\
\leq\rho_\star\sum_{j=0}^{N-1}(1-\beta)^{\lfloor j/N_{\mathrm d}\rfloor}
\leq\frac{N_{\mathrm d}\rho_\star}{\beta},
\end{multline}
where the last inequality follows from $\sum_{j=0}^{N-1}(1-\beta)^{\lfloor j/N_{\mathrm d}\rfloor}\leq N_{\mathrm d}\sum_{k=0}^\infty(1-\beta)^k=N_{\mathrm d}/\beta$. Substituting \eqref{eq:bias_bound_finitehorizon_stationary} and \eqref{eq:relative_bound_finitehorizon_stationary} into \eqref{eq:diff_finitehorizon_stationary_nobounds} proves \eqref{eq:bounded_regret}.

We finally prove (iii). For feasible $\lambda$, $\widetilde\lambda$, and $\tau\in[0,1]$, let $\lambda_\tau:=\tau\lambda+(1-\tau)\widetilde\lambda$. Order concavity gives $
\Phi_\theta(\lambda_\tau)\geq_{\mathcal X^\ast}\tau\Phi_\theta(\lambda)+(1-\tau)\Phi_\theta(\widetilde\lambda)
\geq_{\mathcal X^\ast}\lambda_\tau$.
Hence, the feasible set is convex. Every feasible $\lambda$ satisfies $\Phi_\theta^N(\lambda)\geq_{\mathcal X^\ast}\lambda$ by order preservation. Letting $N\to\infty$ and using (i) and closedness of $\mathcal X^\ast$ gives $\lambda_\star\geq_{\mathcal X^\ast}\lambda$. Since $\lambda_\star$ is feasible and $\xi\in\operatorname{int}\mathcal X$, \eqref{eq:cone_interior_norm} gives uniqueness of the maximizer, as in the proof of Corollary~\ref{cor:supersolution_program}.
\end{proof}
For $\theta\leq0$, \eqref{eq:global_convergence} and \eqref{eq:subsolution_program} in Theorem~\ref{thm:detectability_consequences} provide a geometric convergence bound and a convex characterization of $\lambda_\star$, respectively. In particular, the relative values $J_{N,\lambda_0}(x)-J_{N,\lambda_0}(0)$ converge to $\lambda_\star^\top x$ for any $\lambda_0\in\mathcal X^\ast$.
\begin{remark}
\label{rem:related_work}
Related geometric estimates for deterministic total-cost value iteration are established in \cite[Proposition~2]{Rantzer_2006} which is used in \cite[proof of Theorem~1]{Pates_Rantzer_2024} to prove $\lambda^\top x\leq J_\star(x)$ for feasible $\lambda$ of the convex program for deterministic control on proper cones. Geometric convergence for positive concave maps on the nonnegative orthant is established in \cite[Proposition~3]{Piotrowski_Cavalcante_2022}. For risk-neutral positive semilinear systems, \cite[Propositions~1,3,7--8]{semilineardynamicprogramminganalysis} establishes uniqueness, convergence, and a convex characterization for discounted and undiscounted total-cost problems, using a policy class with linear closed-loop dynamics and costs. Together, \eqref{eq:global_convergence} and \eqref{eq:subsolution_program} extend this combination of geometric convergence and convex characterization to the risk-sensitive average-cost setting on general proper cones, under the hypotheses of Theorem~\ref{thm:detectability_consequences}.
\end{remark}


\begin{remark}
\label{rem:finite_horizon_regret}
For $\theta\leq0$, statement~${\rm (ii)}$ of Theorem~\ref{thm:detectability_consequences} strengthens Remark~\ref{rem:stationary_finite_horizon} by establishing through \eqref{eq:bounded_regret} that $0\leq\sup_{N\geq1}\Delta_N(x)\leq\frac{N_{\mathrm d}\rho_\star}{\beta}+\lambda_\star^\top x<\infty$,
for every fixed $x\in\mathcal X$. Consequently, $\Delta_N(x)/N=O(1/N)$. The state-dependent term in \eqref{eq:bounded_regret} decreases geometrically, although $\Delta_N(x)$ itself need not converge to zero.
\end{remark}

\begin{remark}
\label{rem:risk_averse_nonuniqueness}
For $\theta>0$, order concavity need not hold, and uniqueness may fail even under Assumption~\ref{ass:finite_step_detectability}. In Example~\ref{ex:poisson_assumption} with $n_x=1$, $\theta=1$, $C=1/2$, and $s=1/10$, observability holds with $N_{\mathrm d}=1$ by Remark~\ref{rem:detectability_example1}, even with $\mathcal R_{-|\theta|}$ replaced by $\mathcal R_\theta$ in \eqref{eq:finite_step_detectability}. However, $\Phi_1(\lambda)=1/10+(e^\lambda-1)/2$ has two fixed points, $\lambda_\star\approx0.2281$ and $\lambda_\dagger\approx1.0958$. 
\end{remark}

\subsection{State Bounds}
\label{sec:state_bounds}
For $\theta\geq0$, the following proposition uses Assumption~\ref{ass:finite_step_detectability} to bound $\mathbb E[\|\rvX_t\|]$ uniformly in $t$ under any stationary policy satisfying \eqref{eq:finite_step_stabilizability}. For $\theta>0$, it also gives $\sup_{t\geq0}\mathbb E[e^{(\theta/2)\bar\varepsilon\|\rvX_t\|}]<\infty$ for some $\bar\varepsilon>0$.

\begin{proposition}
\label{prop:finite_step_state_bounds}
Let $\theta\geq0$, and suppose that Assumption~\ref{ass:finite_step_detectability} holds. Let $\bar\lambda\in\mathcal D_\theta$, $\bar\rho\geq0$, and a stationary policy $\bar\pi\in\Pi$ satisfy \eqref{eq:finite_step_stabilizability} for every $N\geq1$ and $x\in\mathcal X$. Then $\bar\lambda\in\operatorname{int}\mathcal X^\ast$, and there exist $\varepsilon\in(0,1]$ and $\bar\varepsilon>0$ such that, for every initial state $\rvX_0=x\in\mathcal X$,
\begin{subequations}
\label{eq:uniform_state_bound}
\begin{equation}
\label{eq:uniform_mean_state_bound}\sup_{t\geq0}\mathbb E[\|\rvX_t\|]\leq\frac{\max\left\{\bar\lambda^\top x,\frac{N_{\mathrm d}\bar\rho}{\varepsilon}\right\}+(N_{\mathrm d}-1)\bar\rho}{\bar\varepsilon}.
\end{equation}
For $\theta>0$, we also have
\begin{multline}
\label{eq:uniform_exponential_state_bound}\sup_{t\geq0}\mathbb E\left[e^{(\theta/2)\bar\varepsilon\|\rvX_t\|}\right]\\ \leq\exp\left[\frac{\theta}{2}\left(\max\left\{\bar\lambda^\top x,\frac{N_{\mathrm d}\bar\rho}{\varepsilon}\right\}+(N_{\mathrm d}-1)\bar\rho\right)\right].
\end{multline}
\end{subequations}
In particular, when Assumption~\ref{ass:rs_continuity} and the equivalent conditions in Theorem~\ref{thm:average_cost} hold, \eqref{eq:uniform_state_bound} applies to the policy in \eqref{eq:opt_sol_star}, with $\bar\lambda=\lambda_\star$ and $\bar\rho=\rho_\star$.
\end{proposition}
\begin{proof}

Since $0\leq S_{N,\bar\pi}^{0}(x)$ and $0\leq\bar\lambda^\top\rvX_N$, monotonicity and \eqref{eq:finite_step_stabilizability} give $0\leq\max\left\{\mathcal R_\theta\left(S_{N,\bar\pi}^{0}(x)\right),\mathcal R_\theta\left(\bar\lambda^\top\rvX_N\right)\right\}\leq\mathcal R_\theta\left(S_{N,\bar\pi}^{\bar\lambda}(x)\right)\leq N\bar\rho+\bar\lambda^\top x<\infty$. Assumption~\ref{ass:finite_step_detectability}, Jensen's inequality, and monotonicity give
\begin{multline*}
\eta^\top x\leq\mathcal R_{-\theta}\left(S_{N_{\mathrm d},\bar\pi}^{0}(x)\right)\leq\mathcal R_\theta\left(S_{N_{\mathrm d},\bar\pi}^{0}(x)\right)\\ \leq\mathcal R_\theta\left(S_{N_{\mathrm d},\bar\pi}^{\bar\lambda}(x)\right)\leq\bar\lambda^\top x+N_{\mathrm d}\bar\rho.
\end{multline*}
Since this holds for every $x\in\mathcal X$ and $\eta\in\operatorname{int}\mathcal X^\ast$, \eqref{eq:rep} and \eqref{eq:cone_interior_addition} give
\begin{equation}
\label{eq:detectability_coefficient_domination}\bar\lambda\geq_{\mathcal X^\ast}\eta,~\bar\lambda\in\operatorname{int}\mathcal X^\ast.
\end{equation}
Moreover, since $\eta\in\operatorname{int}\mathcal X^\ast$, \eqref{eq:cone_interior_comparison} gives $\varepsilon\in(0,1]$ such that $\eta\geq_{\mathcal X^\ast}\varepsilon\bar\lambda$. Hence,
\begin{equation}
\label{eq:functional_domination}\eta^\top x\geq\varepsilon\bar\lambda^\top x,~x\in\mathcal X.
\end{equation}

For $\theta>0$, \eqref{eq:Nstage_cost_terminal} and the Cauchy-Schwarz inequality give
\begin{multline*}
\mathbb E\left[e^{(\theta/2)\bar\lambda^\top\rvX_{N_{\mathrm d}}}\right]=\mathbb E\left[e^{(\theta/2)S_{N_{\mathrm d},\bar\pi}^{\bar\lambda}(x)}e^{-(\theta/2)S_{N_{\mathrm d},\bar\pi}^{0}(x)}\right]\\ \leq\left(\mathbb E\left[e^{\theta S_{N_{\mathrm d},\bar\pi}^{\bar\lambda}(x)}\right]\right)^{1/2}\left(\mathbb E\left[e^{-\theta S_{N_{\mathrm d},\bar\pi}^{0}(x)}\right]\right)^{1/2}.
\end{multline*}
Taking logarithms and multiplying by $2/\theta$ yields
\begin{equation}
\label{eq:finite_step_risk_separation}\mathcal R_{\theta/2}\left(\bar\lambda^\top\rvX_{N_{\mathrm d}}\right)\leq\mathcal R_\theta\left(S_{N_{\mathrm d},\bar\pi}^{\bar\lambda}(x)\right)-\mathcal R_{-\theta}\left(S_{N_{\mathrm d},\bar\pi}^{0}(x)\right).
\end{equation}
For $\theta=0$, this holds with equality by linearity of expectation. Combining \eqref{eq:finite_step_risk_separation} with \eqref{eq:finite_step_stabilizability}, \eqref{eq:finite_step_detectability}, and \eqref{eq:functional_domination} gives
\begin{equation}
\label{eq:finite_step_state_nonnegative}\mathcal R_{\theta/2}\left(\bar\lambda^\top\rvX_{N_{\mathrm
d}}\right)\leq\bar\lambda^\top x-\eta^\top x+N_{\mathrm d}\bar\rho\leq(1-\varepsilon)\bar\lambda^\top x+N_{\mathrm d}\bar\rho.
\end{equation}

The convexity argument in the proof of Proposition~\ref{prop:convexthetaR}, applied to $\tau\bar\lambda^\top\rvX_t+(1-\tau)0$, together with $\mathcal R_\theta(0)=0$, gives
\begin{equation}
\label{eq:nonnegative_risk_scaling}\mathcal R_\theta\left(\tau\bar\lambda^\top\rvX_t\right)\leq\tau\mathcal R_\theta\left(\bar\lambda^\top\rvX_t\right),
\end{equation}
for $\theta>0$, $\tau\in[0,1]$, and $t\geq0$. For $\theta=0$, the same inequality holds with equality by linearity of expectation.

For every $k\geq0$, stationarity of $\bar\pi$ and \eqref{eq:general_dynamics} imply that, conditional on $\rvX_{kN_{\mathrm d}}$, the next $N_{\mathrm d}$ stages have the law induced by $\rvX_{kN_{\mathrm d}}$, $\bar\pi$, and $\rcd$. Therefore, \eqref{eq:finite_step_state_nonnegative}, applied with initial state $\rvX_{kN_{\mathrm d}}$, gives
$\mathcal R_{\theta/2}\left(\bar\lambda^\top\rvX_{(k+1)N_{\mathrm d}}\mid\rvX_{kN_{\mathrm d}}\right)\leq(1-\varepsilon)\bar\lambda^\top\rvX_{kN_{\mathrm d}}+N_{\mathrm d}\bar\rho$.
Using \eqref{eq:recursive_R}, monotonicity, translation invariance, and \eqref{eq:nonnegative_risk_scaling} at risk parameter $\theta/2$, we obtain
\begin{multline*}
\mathcal R_{\theta/2}\left(\bar\lambda^\top\rvX_{(k+1)N_{\mathrm d}}\right)=\mathcal R_{\theta/2}\left(\mathcal R_{\theta/2}\left(\bar\lambda^\top\rvX_{(k+1)N_{\mathrm d}}\mid\rvX_{kN_{\mathrm d}}\right)\right)\\ \leq\mathcal R_{\theta/2}\left((1-\varepsilon)\bar\lambda^\top\rvX_{kN_{\mathrm d}}\right)+N_{\mathrm d}\bar\rho\\ \leq(1-\varepsilon)\mathcal R_{\theta/2}\left(\bar\lambda^\top\rvX_{kN_{\mathrm d}}\right)+N_{\mathrm d}\bar\rho.
\end{multline*}
Since $\rvX_0=x$, iterating this inequality gives, for every $k\geq1$,
\begin{multline*}
\mathcal R_{\theta/2}\left(\bar\lambda^\top\rvX_{kN_{\mathrm d}}\right)\leq(1-\varepsilon)^k\bar\lambda^\top x+N_{\mathrm d}\bar\rho\sum_{i=0}^{k-1}(1-\varepsilon)^i\\ =(1-\varepsilon)^k\bar\lambda^\top x+\left(1-(1-\varepsilon)^k\right)\frac{N_{\mathrm d}\bar\rho}{\varepsilon} \leq\max\left\{\bar\lambda^\top x,\frac{N_{\mathrm d}\bar\rho}{\varepsilon}\right\}.
\end{multline*}
The last inequality follows from $0\leq(1-\varepsilon)^k\leq1$, and the same bound holds for $k=0$.

It remains to bound the state at stages between successive multiples of $N_{\mathrm d}$. By Jensen's inequality, nonnegativity of the stage cost, and \eqref{eq:finite_step_stabilizability} with $N=1$,
$\mathcal R_{\theta/2}\left(\bar\lambda^\top\rvX_{t+1}\mid\rvX_t\right)\leq\mathcal R_\theta\left(\bar\lambda^\top\rvX_{t+1}\mid\rvX_t\right)\leq\bar\rho+\bar\lambda^\top\rvX_t$.
Applying \eqref{eq:recursive_R}, monotonicity, and translation invariance successively gives, for $k\geq0$ and $0\leq j<N_{\mathrm d}$,
\begin{multline*}
\mathcal R_{\theta/2}\left(\bar\lambda^\top\rvX_{kN_{\mathrm d}+j}\right)\leq\mathcal R_{\theta/2}\left(\bar\lambda^\top\rvX_{kN_{\mathrm d}}\right)+j\bar\rho\\ \leq\max\left\{\bar\lambda^\top x,\frac{N_{\mathrm d}\bar\rho}{\varepsilon}\right\}+(N_{\mathrm d}-1)\bar\rho.
\end{multline*}
Since every $t\geq0$ can be written as $t=kN_{\mathrm d}+j$ with $k\geq0$ and $0\leq j<N_{\mathrm d}$, this bound holds for every $t\geq0$.

Since $\bar\lambda\in\operatorname{int}\mathcal X^\ast$, \eqref{eq:cone_interior_norm} gives $\bar\varepsilon>0$ such that $\bar\lambda^\top x\geq\bar\varepsilon\|x\|$ for every $x\in\mathcal X$. Jensen's inequality gives $\bar\varepsilon\mathbb E[\|\rvX_t\|]\leq\mathbb E[\bar\lambda^\top\rvX_t]\leq\mathcal R_{\theta/2}(\bar\lambda^\top\rvX_t)$. For $\theta>0$, we also have $\mathbb E[e^{(\theta/2)\bar\varepsilon\|\rvX_t\|}]\leq\mathbb E[e^{(\theta/2)\bar\lambda^\top\rvX_t}]=\exp((\theta/2)\mathcal R_{\theta/2}(\bar\lambda^\top\rvX_t))$. Combining these inequalities with the preceding bound proves \eqref{eq:uniform_state_bound}.

Finally, if Assumption~\ref{ass:rs_continuity} and the equivalent conditions in Theorem~\ref{thm:average_cost} hold, then \eqref{eq:terminal_identity} gives \eqref{eq:finite_step_stabilizability} with equality for $\bar\lambda=\lambda_\star$, $\bar\rho=\rho_\star$, and $\bar\pi=\pi_\star$. Hence, the bounds apply to the policy in \eqref{eq:opt_sol_star}.
\end{proof}

\begin{remark}
\label{rem:conditional_cost_observability}
For $\theta>0$, the conclusions of Proposition~\ref{prop:finite_step_state_bounds} hold with $\theta/2$ replaced by $\theta$ in \eqref{eq:uniform_exponential_state_bound} if Assumption~\ref{ass:finite_step_detectability} is replaced by
\begin{equation}
\label{eq:conditional_cost_observability}
\mathcal R_0\left(S_{N_{\mathrm d},\pi_{N_{\mathrm d}}}^{0}(x)\mid\rvX_{N_{\mathrm d}}\right)\geq\eta^\top x
\end{equation}
for some $N_{\mathrm d}\geq1$, $\eta\in\operatorname{int}\mathcal X^\ast$, and every $x\in\mathcal X$, $\pi_{N_{\mathrm d}}\in\Pi_{N_{\mathrm d}}$. By Jensen's inequality
\begin{equation}
\label{eq:finite_step_detectability_positive}
\mathcal R_\theta\left(S_{N_{\mathrm d},\pi_{N_{\mathrm d}}}^{\lambda}(x)\right)\geq\eta^\top x+\mathcal R_\theta\left(\lambda^\top\rvX_{N_{\mathrm d}}\right).
\end{equation}
Equation~\eqref{eq:conditional_cost_observability} gives \eqref{eq:detectability_coefficient_domination}, while \eqref{eq:finite_step_detectability_positive} gives \eqref{eq:finite_step_state_nonnegative} with $\mathcal R_\theta$ in place of $\mathcal R_{\theta/2}$. The remaining argument using \eqref{eq:nonnegative_risk_scaling} then gives the claimed bound.
\end{remark}

For $\theta<0$, the inequality in \eqref{eq:nonnegative_risk_scaling} is reversed, so the argument used to obtain \eqref{eq:uniform_state_bound} does not apply.

\section{Connection to Deterministic Infsup Control}
\label{sec:infsup_connection}
In this section, we show that when the process noise enters the dynamics and the stage cost affienly
\begin{subequations}
\label{eq:infsup_model}
\begin{align}
&f(x,u,w)=f_0(x,u)+F(x,u)w,\\
&\ell(x,u,w)=\ell_0(x,u)+\ell_w^\top(x,u)w,
\end{align}
\end{subequations}
where $F(x,u)\in\mathbb R^{n_x\times n_w}$ and $\ell_w(x,u)\in\mathbb R^{n_w}$, then the risk-averse problem ($\theta>0$) is equivilant to a deterministic infsup problem with affine value functions.\footnote{We will focus on the risk-averse problem. The risk-seeeking problem can also be shown to be equivilant to an infinf problem with similar arguments.} 
Note that the models in Examples~\ref{ex:polyhedral_ex1intro} and~\ref{ex:poisson_assumption} have the form in \eqref{eq:infsup_model}. For $J(x)=\lambda^\top x$ with $\lambda\in\mathcal X^\ast$, substituting \eqref{eq:infsup_model} into \eqref{eq:risk_functional} gives
\begin{multline}
\label{eq:infsup_split}\mathcal R_\theta\left(\mathcal H(x,u,\rvW,\lambda)\mid x,u\right)=\ell_0(x,u)+\lambda^\top f_0(x,u)\\+\frac{1}{\theta}\Psi_{x,u}\big(\theta v_\lambda(x,u)\big),
\end{multline}
where $v_\lambda(x,u):=\ell_w(x,u)+F^\top(x,u)\lambda$ and $\Psi_{x,u}(v):=\mathcal R_1(v^\top\rvW\mid x,u)$ is the conditional log moment-generating function. Define its convex conjugate by
\begin{equation}
\label{eq:rate_function}\mathcal I_{x,u}(d):=\sup_{v\in\mathbb R^{n_w}}\left\{v^\top d-\Psi_{x,u}(v)\right\},
\end{equation}
and let $\mathcal D_{\mathcal I}(x,u):=\{d\in\mathbb R^{n_w}:\mathcal I_{x,u}(d)<\infty\}$. The function in \eqref{eq:rate_function} is the Cram\'er rate function of $\rcd(\cdot\mid x,u)$. It is nonnegative and vanishes at $\mathbb E[\rvW\mid x,u]$ whenever the conditional mean is finite. The function $\Psi_{x,u}$ is convex and lower semicontinuous. It is also proper, since $\Psi_{x,u}(0)=0$ and $\Psi_{x,u}(v)>-\infty$ for every $v$. Hence, by the Fenchel-Moreau theorem, we have
\begin{equation}
\label{eq:infsup_duality}\frac{1}{\theta}\Psi_{x,u}(\theta v)=\sup_{d\in\mathcal D_{\mathcal I}(x,u)}\left\{v^\top d-\frac{1}{\theta}\mathcal I_{x,u}(d)\right\}.
\end{equation}
Substituting \eqref{eq:infsup_duality} into \eqref{eq:infsup_split}, using \eqref{eq:infsup_model}, and minimizing over $u\in\mathcal A(x)$ as in \eqref{eq:Bellman_operator}, we obtain
\begin{multline}
\label{eq:infsup_bellman}(T_\theta J)(x)=\inf_{u\in\mathcal A(x)}\sup_{d\in\mathcal D_{\mathcal I}(x,u)}\Big\{\ell(x,u,d)-\frac{1}{\theta}\mathcal I_{x,u}(d)\\+\lambda^\top f(x,u,d)\Big\}.
\end{multline}
The right-hand side of \eqref{eq:infsup_bellman} is the Bellman operator, applied to $J(x)=\lambda^\top x$, of the following deterministic discrete-time infsup problem \cite{Bertsekas_2022_ADP}
\begin{multline}
\label{eq:infsup_problem}
\inf_{\pi_N\in\Pi_N}\sup_{d_0,\ldots,d_{N-1}}\Bigg\{\sum_{t=0}^{N-1}\Big(\ell(x_t,u_t,d_t)\\
-\frac{1}{\theta}\mathcal I_{x_t,u_t}(d_t)\Big)+\lambda_0^\top x_N\Bigg\},
\end{multline}
subject to $x_0=x$, $x_{t+1}=f(x_t,u_t,d_t)$, and $u_t=a_t(x_t)$, with the supremum taken over disturbances satisfying $d_t\in\mathcal D_{\mathcal I}(x_t,u_t)$ for all $t\geq0$. The term $\theta^{-1}\mathcal I_{x_t,u_t}(d_t)$ in \eqref{eq:infsup_problem} penalizes the disturbance.
We now state the following proposition
\begin{proposition}
\label{prop:infsup_admissibility}
Suppose that \eqref{eq:infsup_model}. Then for every $x\in\mathcal X$, $u\in\mathcal A(x)$, and $d\in\mathcal D_{\mathcal I}(x,u)$, we have $f(x,u,d)\in\mathcal X$, and $\ell(x,u,d)\geq0$.
\end{proposition}
\begin{proof}
Fix $x\in\mathcal X$ and $u\in\mathcal A(x)$, and let $C_{x,u}:=\{d\in\mathbb R^{n_w}:f(x,u,d)\in\mathcal X,\ \ell(x,u,d)\geq0\}$. Since $f$ and $\ell$ are affine in $d$ by \eqref{eq:infsup_model}, and $\mathcal X$ is closed and convex, the set $C_{x,u}$ is closed and convex. By admissibility and Assumption~\ref{ass:basic}, $\rcd(C_{x,u}\cap\mathcal W\mid x,u)=1$. Since the support of a probability measure is the smallest closed set of probability one, $C_{x,u}$ contains the support of $\rcd(\cdot\mid x,u)$, viewed as a probability measure on $\mathbb R^{n_w}$, and hence its closed convex hull, which contains $\mathcal D_{\mathcal I}(x,u)$ by \cite[Proposition~1(a)]{Vysotsky_2021} and \eqref{eq:rate_function}.
\end{proof}
Note that by \eqref{eq:rs_effective_domain} and \eqref{eq:infsup_bellman}, $\mathcal D_\theta$ is also the set of $\lambda\in\mathcal X^\ast$ for which the deterministic infsup operator in \eqref{eq:infsup_bellman} applied to $J(x)=\lambda^\top x$ is finite for every $x\in\mathcal X$. The following proposition shows that \eqref{eq:infsup_problem} has the same finite-horizon value and optimal policy as the corresponding risk-sensitive problem.
\begin{proposition}
\label{prop:infsup_finite_horizon}
Suppose that \eqref{eq:infsup_model} holds and that the recurssion in \eqref{eq:affine_recursion}, initialized with $c_0=0$ and $\lambda_0\in\mathcal{D}_\theta$, is well defined up to horizon $N$. Then the value of \eqref{eq:infsup_problem} is $c_N+\lambda_N^\top x$ for every $x\in\mathcal X$, and its infimum is attained by the policy in \eqref{eq:finite_horizon_optimal_policy}.
\end{proposition}
\begin{proof}
For $J(x)=c+\lambda^\top x$ with $\lambda\in\mathcal D_\theta$, \eqref{eq:infsup_split}, \eqref{eq:infsup_duality}, and \eqref{eq:bellman_affinity} show that both Bellman operators give $c+\Gamma_\theta(\lambda)+\Phi_\theta(\lambda)^\top x$, with the same minimizing decision rule in \eqref{eq:opt_sol}. Proposition~\ref{prop:infsup_admissibility} ensures that $f(x,u,d)\in\mathcal X$ for $u\in\mathcal A(x)$ and $d\in\mathcal D_{\mathcal I}(x,u)$, so the induction in the proof of Proposition~\ref{prop:finite_horizon_affine} applies.
\end{proof}
We now provide a proposition for the corresponding time-average problem.
\begin{proposition}
\label{prop:infsup_average_cost}
Suppose that \eqref{eq:infsup_model} holds together with Assumption~\ref{ass:rs_continuity} and the equivalent conditions of Theorem~\ref{thm:average_cost}. Then, for every $x\in\mathcal X$,
\begin{multline}
\label{eq:infsup_average_cost}
\inf_{\pi\in\Pi}\limsup_{N\to\infty}\frac{1}{N}\sup_{d_0,\ldots,d_{N-1}}
\sum_{t=0}^{N-1}\Big(\ell(x_t,u_t,d_t)\\
-\frac{1}{\theta}\mathcal I_{x_t,u_t}(d_t)\Big)=\rho_\star,
\end{multline}
where $x_0=x$, $x_{t+1}=f(x_t,u_t,d_t)$, $u_t=a_t(x_t)$, and the supremum is over disturbances satisfying $d_t\in\mathcal D_{\mathcal I}(x_t,u_t)$ for all $t\geq0$. The infimum is attained by the stationary policy in \eqref{eq:opt_sol_star}.
\end{proposition}
\begin{proof}
Fix $x\in\mathcal X$ and let $\pi\in\Pi$ be arbitrary. By Proposition~\ref{prop:infsup_finite_horizon} with $\lambda_0=0$,
\begin{equation}
\label{eq:finite-horizon-sup-lowerbound}
\sup_{d_0,\ldots,d_{N-1}}\sum_{t=0}^{N-1}\Big(\ell(x_t,u_t,d_t)
-\frac{1}{\theta}\mathcal I_{x_t,u_t}(d_t)\Big)
\geq c_N+\lambda_N^\top x.
\end{equation}
Dividing \eqref{eq:finite-horizon-sup-lowerbound} by $N$, taking the upper limit, and using $\lambda_N\to\lambda_\star$ and $c_N/N\to\rho_\star$ from Theorem~\ref{thm:average_cost} and its proof proves the lower bound in \eqref{eq:infsup_average_cost}.

For the reverse inequality, let $\pi_\star$ be the stationary policy in \eqref{eq:opt_sol_star}. By \eqref{eq:fixed_point_bellman} and \eqref{eq:infsup_bellman}, every admissible disturbance sequence under $\pi_\star$ satisfies
$\ell(x_t,u_t,d_t)-\frac{1}{\theta}\mathcal I_{x_t,u_t}(d_t)
+\lambda_\star^\top x_{t+1}
\leq\rho_\star+\lambda_\star^\top x_t$.
Summing over $t\in\{0,\ldots,N-1\}$ gives
\begin{equation}
\label{eq:infsup_dissipation}
\sum_{t=0}^{N-1}\Big(\ell(x_t,u_t,d_t)-\frac{1}{\theta}\mathcal I_{x_t,u_t}(d_t)\Big)
+\lambda_\star^\top x_N\leq N\rho_\star+\lambda_\star^\top x.
\end{equation}
By Proposition~\ref{prop:infsup_admissibility}, $x_N\in\mathcal X$, so $\lambda_\star^\top x_N\geq0$. Taking the supremum over disturbances in \eqref{eq:infsup_dissipation}, dividing by $N$, and taking the upper limit therefore shows that the time-average cost under $\pi_\star$ is at most $\rho_\star$. Together with the lower bound, this proves \eqref{eq:infsup_average_cost} and optimality of $\pi_\star$.
\end{proof}

\begin{remark}
\label{rem:infsup_total_cost}
Under \eqref{eq:infsup_model} and the assumptions of Proposition~\ref{prop:stochastic_total_cost}, we have $c_N=\rho_\star=0$. Arguments similar to its proof together with \eqref{eq:finite-horizon-sup-lowerbound} and \eqref{eq:infsup_dissipation} give
\begin{equation}
\label{eq:infsup_total_cost}\inf_{\pi\in\Pi}\sup_{N\geq1}\sup_{d_0,\ldots,d_{N-1}}\sum_{t=0}^{N-1}\Big(\ell(x_t,u_t,d_t)-\frac{1}{\theta}\mathcal I_{x_t,u_t}(d_t)\Big)=\lambda_\star^\top x,
\end{equation}
for every $x\in\mathcal{X}$ and $d_t\in\mathcal{D}_I$ for all $t\geq 0$. Here, the infimum is attained by the stationary policy in \eqref{eq:opt_sol_star}.
\end{remark}

For Example~\ref{ex:polyhedral_ex1intro}, the IG distribution in \eqref{eq:polyhedral_ig} gives $\mathcal I_{x,u}(d)=\kappa^\top d-2\kappa^\top\mu_w(x,u)+\sum_{i=1}^{n_d}\kappa^{(i)}(\mu_w^{(i)}(x,u))^2/d^{(i)}$ for $d>0$, and $\mathcal I_{x,u}(d)=\infty$ otherwise. The corresponding infsup problem \eqref{eq:infsup_problem} has dynamics $x_{t+1}=Ax_t+Bu_t+Fd_t$ with $u_t\in\mathcal U(x_t)$ as in Example~\ref{ex:polyhedral_ex1intro}, disturbances $d_t>0$, and stage cost $s^\top x+r^\top u-\theta^{-1}\mathcal I_{x,u}(d)$.
For fixed $(x,u)$, we have that the supremum in \eqref{eq:infsup_bellman} is finite if and only if $\theta(F^\top\lambda)\leq\kappa$, which  recovers $\mathcal D_\theta$ in Example~\ref{ex:polyhedral_ex1intro}. The same finiteness condition holds with the linear penalty $\gamma^\top d$ for $d\geq0$, where $\gamma=\kappa/\theta$, as in \cite{Gurpegui_Tegling_Rantzer_2025_Unbounded}. Moreover, for fixed $d>0$, $\theta^{-1}\mathcal I_{x,u}(d)\to\gamma^\top d$ as $\mu_w(x,u)\to0$.

\section{Additional Examples}
\label{sec:examples}
In this section, we present two additional examples to the one discussed already. The first one is for PSD cones, and the second one is a numerical example.

\subsection{PSD Cone and LEQG}
\label{sec:PSDexample}
Let $\mathbb S^{n_x}$ denote the space of symmetric $n_x\times n_x$ matrices equipped with the inner product $\langle X,Y\rangle:=\operatorname{tr}(XY)$, and take $\mathcal X=\mathbb S_{\geq}^{n_x}$ to be the cone of PSD matrices.\footnote{As in \cite[Sec.~3.1]{Pates_Rantzer_2024}, we use matrix notation rather than explicitly identifying $\mathbb S^{n_x}$ with $\mathbb R^{n_x(n_x+1)/2}$. The results for proper cones apply equally under the trace inner product.}

Since $\mathbb S_{\geq}^{n_x}$ is self-dual, one has $\mathcal X^\ast=\mathbb S_{\geq}^{n_x}$. In other words, for $\Lambda\in\mathcal X^\ast$, we have $\operatorname{tr}(\Lambda X)\geq0$ for any $X\in\mathcal X$. Let
\begin{multline*}
\mathbb U_{n_x,n_u}=\left\{U=\begin{bmatrix}0&U_1\\U_1^\top&U_2\end{bmatrix}:\right.\\\left.U_1\in\mathbb R^{n_x\times n_u},~U_2\in\mathbb S^{n_u}\right\},
\end{multline*}
and for every $X\in\mathcal X$, define
\begin{multline}
\label{eq:psd_input_set}\mathcal U(X):=\Biggl\{U\in\mathbb U_{n_x,n_u}:\begin{bmatrix}X&0\\0&0\end{bmatrix}+U\in\mathbb S_{\geq}^{n_x+n_u}\Biggr\}.
\end{multline}

Note that the set $\mathcal P=\left\{(X,U):X\in\mathcal X,~U\in\mathcal U(X)\right\}$, equipped with the inner product $\langle(X,U),(H,M)\rangle:=\operatorname{tr}(XH)+\operatorname{tr}(UM)$ on $\mathbb S^{n_x}\times\mathbb U_{n_x,n_u}$, is a proper cone and it is self-dual, i.e., $\mathcal P^\ast=\mathcal P$.

We consider dynamics given by
\begin{equation}
\label{eq:psd_stochastic_dynamics}\mrvX_{t+1}=\mrvW_t,
\end{equation}
where $\mrvW_t\in\mathbb S_{\geq}^{n_x}$. Hence, we have that $\mathcal A(X)=\mathcal U(X)$. Furthermore, we choose the stage cost to be
\begin{equation}
\label{eq:psd_stage_cost}\ell(X,U)=\operatorname{tr}(QX)+\operatorname{tr}(RU),
\end{equation}
where $(Q,R)\in\mathcal P^\ast$. Additionally, for $(Q,R)\in\operatorname{int}\mathcal P^\ast$, Assumption~\ref{ass:finite_step_detectability} holds with $N_{\mathrm d}=1$ (see Remark~\ref{rem:intpropercone_detect}).

Now let $\rcd$ be any regular conditional distribution of $\mrvW_t$ given $(\mrvX_t,\mrvU_t)$ such that, for any $\Lambda\in\mathcal D_\theta$, we have
\begin{multline}
\label{eq:psd_w2condition}\mathcal R_\theta\left(\operatorname{tr}(\Lambda\mrvW_t)\mid X,U\right)=\\\nu_\theta(\Lambda)+\operatorname{tr}(N_{\theta,x}(\Lambda)X)+\operatorname{tr}(N_{\theta,u}(\Lambda)U)
\end{multline}
for all $(X,U)\in\mathcal P$, where $\nu_\theta(\Lambda)\in\mathbb R$, $N_{\theta,x}(\Lambda)\in\mathbb S^{n_x}$, and $N_{\theta,u}(\Lambda)\in\mathbb U_{n_x,n_u}$. Since $\operatorname{tr}(\Lambda\mrvW_t)\geq0$ for any $\Lambda\in\mathcal D_\theta$, we have that the left-hand side of \eqref{eq:psd_w2condition} is nonnegative, which implies according to \eqref{eq:rep} that $\nu_\theta(\Lambda)\geq0$ and $(N_{\theta,x}(\Lambda),N_{\theta,u}(\Lambda))\in\mathcal P^\ast=\mathcal P$ for any $\Lambda\in\mathcal D_\theta$.

For this example, we have that, for each $\Lambda\in\mathcal D_\theta$,
\begin{multline}
\label{eq:psd_RH}\mathcal R_\theta\left(\mathcal H(X,U,\mrvW,\Lambda)\mid X,U\right)=\\\nu_\theta(\Lambda)+\operatorname{tr}(H_\theta(\Lambda)X)+\operatorname{tr}(M_\theta(\Lambda)U),
\end{multline}
where $H_\theta(\Lambda)=Q+N_{\theta,x}(\Lambda)$ and $M_\theta(\Lambda)=R+N_{\theta,u}(\Lambda)$. Additionally, since $(Q,R)\in\mathcal P^\ast=\mathcal P$ and $(N_{\theta,x}(\Lambda),N_{\theta,u}(\Lambda))\in\mathcal P$, we have that $(H_\theta(\Lambda),M_\theta(\Lambda))\in\mathcal P$. Writing
$$M_\theta(\Lambda)=\begin{bmatrix}0&M_{\theta,1}(\Lambda)\\M_{\theta,1}^\top(\Lambda)&M_{\theta,2}(\Lambda)\end{bmatrix},~R=\begin{bmatrix}0&R_1\\R_1^\top&R_2\end{bmatrix},$$
we have in particular that $M_{\theta,2}(\Lambda)\in\mathbb S_{\geq}^{n_u}$ and $\operatorname{range}(M_{\theta,1}^\top(\Lambda))\subseteq\operatorname{range}(M_{\theta,2}(\Lambda))$. Then,
\begin{multline}
\label{eq:minimum_psd}\inf_{U:(X,U)\in\mathcal P}\operatorname{tr}(M_\theta(\Lambda)U)=\\-\operatorname{tr}\left(M_{\theta,1}(\Lambda)M_{\theta,2}^\dagger(\Lambda)M_{\theta,1}^\top(\Lambda)X\right),
\end{multline}
and one minimizer is
\begin{multline}
\label{eq:generalminimizer_psd}U=a(X;\Lambda)=\left[\begin{smallmatrix}0&-XM_{\theta,1}M_{\theta,2}^\dagger\\-M_{\theta,2}^\dagger M_{\theta,1}^\top X&M_{\theta,2}^\dagger M_{\theta,1}^\top XM_{\theta,1}M_{\theta,2}^\dagger\end{smallmatrix} \right],
\end{multline}
where $M_{\theta,2}^\dagger$ denotes the pseudoinverse of $M_{\theta,2}$, and the dependence of $M_{\theta,1}$ and $M_{\theta,2}$ on $\Lambda$ is omitted for brevity. Note that if we choose $R_2\in\mathbb S_{>}^{n_u}$ (automatically satisfied for $(Q,R)\in\operatorname{int}\mathcal P$ by \eqref{eq:cone_interior_comparison}), then $M_{\theta,2}\in\mathbb S_{>}^{n_u}$ by \eqref{eq:cone_interior_addition}, so it is invertible and we can replace $M_{\theta,2}^\dagger$ in \eqref{eq:minimum_psd} and \eqref{eq:generalminimizer_psd} with $M_{\theta,2}^{-1}$. Moreover, with $R_2\in\mathbb S_{>}^{n_u}$, the minimizer in \eqref{eq:generalminimizer_psd} is the unique minimizer. With this, we have that Assumption~\ref{ass:affine_bellman} is satisfied with $\Phi_\theta(\Lambda)=H_\theta(\Lambda)-M_{\theta,1}(\Lambda)M_{\theta,2}^\dagger(\Lambda)M_{\theta,1}^\top(\Lambda)$ and $\Gamma_\theta(\Lambda)=\nu_\theta(\Lambda)$.

One example of a process $\mrvW_t$ satisfying \eqref{eq:psd_w2condition} is one where $\mrvW_t\mid(\mrvX_t,\mrvU_t)=(X,U)$ is distributed according to a noncentral Wishart distribution \cite{Anderson_1946,Muirhead_1982}. Specifically, for $(X,U)\in\mathcal P$, let
\begin{equation}
\label{eq:noncentrality_psd}C(X,U):=\begin{bmatrix}A&B\end{bmatrix}\begin{bmatrix}X&U_1\\U_1^\top&U_2\end{bmatrix}\begin{bmatrix}A^\top\\B^\top\end{bmatrix},
\end{equation}
with $A\in\mathbb R^{n_x\times n_x}$ and $B\in\mathbb R^{n_x\times n_u}$. Consider a positive integer $p$ and a measurable factor $\tilde C(X,U)\in\mathbb R^{n_x\times p}$ such that $C(X,U)=\tilde C(X,U)\tilde C^\top(X,U)$. Such a factor can always be chosen for $p\geq\operatorname{rank}(C(X,U))$. Consider $p$ independent Gaussian vectors $\tilde{\mrvW}_{1,t}\sim\mathcal N(0,\Sigma),\ldots,\tilde{\mrvW}_{p,t}\sim\mathcal N(0,\Sigma)$, with $\Sigma\in\mathbb S_{\geq}^{n_x}$, independent of the current state and input and of past disturbances. Then, conditionally on $(\mrvX_t,\mrvU_t)=(X,U)$, the matrix $\mrvW_t=\sum_{i=1}^p(\tilde C_i(X,U)+\tilde{\mrvW}_{i,t})(\tilde C_i(X,U)+\tilde{\mrvW}_{i,t})^\top$, with $\tilde C_i$ denoting the $i$th column of $\tilde C$, has a noncentral Wishart distribution with parameters $p$, $\Sigma$, and noncentrality matrix $C(X,U)$,
$$\mrvW_t\mid(\mrvX_t,\mrvU_t)=(X,U)\sim\mathcal W_{n_x}(p,\Sigma,C(X,U)).$$
The conditional mean is $\mathbb E[\mrvW_t\mid X,U]=p\Sigma+C(X,U)$. Additionally, we have that $\mathcal R_\theta(\operatorname{tr}(\Lambda\mrvW_t)\mid X,U)$ is in the form of \eqref{eq:psd_w2condition} with
\begin{align*}
&\begin{multlined}\nu_\theta(\Lambda)=\begin{cases}-\dfrac{p}{2\theta}\log\det\left(I-2\theta\Lambda^{1/2}\Sigma\Lambda^{1/2}\right),&\theta\neq0,\\p\operatorname{tr}(\Sigma\Lambda),&\theta=0,\end{cases}\end{multlined},\\&N_{\theta,u}(\Lambda)=\begin{bmatrix}0&A^\top\Lambda_{\theta,\Sigma}B\\B^\top\Lambda_{\theta,\Sigma}A&B^\top\Lambda_{\theta,\Sigma}B\end{bmatrix},~N_{\theta,x}(\Lambda)=A^\top\Lambda_{\theta,\Sigma}A
\end{align*}
where
\begin{multline*}
\Lambda_{\theta,\Sigma}:=\begin{cases}\Lambda^{1/2}\left(I-2\theta\Lambda^{1/2}\Sigma\Lambda^{1/2}\right)^{-1}\Lambda^{1/2},&\theta\neq0,\\\Lambda,&\theta=0,\end{cases}
\end{multline*}
provided that $I-2\theta\Lambda^{1/2}\Sigma\Lambda^{1/2}\in\mathbb S_{>}^{n_x}$, and hence
$$\mathcal D_\theta=\{\Lambda\in\mathbb S_{\geq}^{n_x}:I-2\theta\Lambda^{1/2}\Sigma\Lambda^{1/2}\in\mathbb S_{>}^{n_x}\}.$$
Therefore, with the noncentral Wishart choice for $\mrvW_t$, we get that
\begin{subequations}
\begin{align}
&\begin{multlined}\Phi_\theta(\Lambda)=Q+A^\top\Lambda_{\theta,\Sigma}A-\left(A^\top\Lambda_{\theta,\Sigma}B+R_1\right)\\\left(B^\top\Lambda_{\theta,\Sigma}B+R_2\right)^\dagger\left(B^\top\Lambda_{\theta,\Sigma}A+R_1^\top\right),\end{multlined}\label{eq:Wishart_Phi}\\&\Gamma_\theta(\Lambda)=\nu_\theta(\Lambda).\label{eq:Wishart_Gamma}
\end{align}
\end{subequations}
Note, however, that continuity of $\Phi_\theta$ in \eqref{eq:Wishart_Phi} is not automatic when $M_{\theta,2}(\Lambda)=B^\top\Lambda_{\theta,\Sigma}B+R_2$ is singular, since its rank may vary with $\Lambda$. Thus, Assumption~\ref{ass:rs_continuity} must be verified separately. A simple sufficient condition is $R_2\in\mathbb S_{>}^{n_u}$, in which case \eqref{eq:cone_interior_addition} gives $M_{\theta,2}(\Lambda)=B^\top\Lambda_{\theta,\Sigma}B+R_2\in\mathbb S_{>}^{n_u}$ for every $\Lambda\in\mathcal D_\theta$, and $\Phi_\theta$ is continuous.

The equation $\Lambda=\Phi_\theta(\Lambda)$ with $\Phi_\theta$ as in \eqref{eq:Wishart_Phi} is the Riccati equation for LEQG \cite{Jacobson_1973,Whittle_1981,Whittle_1990}. We can further clarify the link with LEQG by following the construction in \cite[Sec.~3.1]{Pates_Rantzer_2024}, where the deterministic LQR problem is recovered from an optimal control problem on the PSD cone. Specifically, consider the following dynamics
\begin{equation}
\label{eq:LQG_dynamics}\rvX_{t+1}=A\rvX_t+B\rvU_t+\rvW_t,
\end{equation}
with $\rvX_t\in\mathbb R^{n_x}$, $\rvU_t\in\mathbb R^{n_u}$, and independent Gaussian disturbances $\rvW_t\sim\mathcal N(0,\Sigma)$, with the cost
\begin{equation}
\label{eq:LQG_Cost}\ell(x,u,w)=x^\top Qx+2x^\top R_1u+u^\top R_2u,
\end{equation}
where $(Q,R)\in\mathcal P^\ast$ and $R_2\in\mathbb S_{>}^{n_u}$. As in the deterministic construction in \cite{Pates_Rantzer_2024}, define the lifted variables $X_t=x_tx_t^\top$, $U_{1,t}=x_tu_t^\top$, and $U_{2,t}=u_tu_t^\top$. For these state-input pairs, \eqref{eq:noncentrality_psd} gives $C(X,U)=(Ax+Bu)(Ax+Bu)^\top$, and hence $\operatorname{rank}(C(X,U))\leq p=1$. Consequently, with $\mrvX_t=\rvX_t\rvX_t^\top$ and $\mrvU_t$ defined through the corresponding lifted inputs, the dynamics in \eqref{eq:LQG_dynamics} can be written in the form of \eqref{eq:psd_stochastic_dynamics} with
$$\mrvW_t\mid(\mrvX_t,\mrvU_t)=(X,U)\sim\mathcal W_{n_x}(1,\Sigma,C(X,U)),$$
and the stage cost in \eqref{eq:LQG_Cost} can be written in the form of the stage cost in \eqref{eq:psd_stage_cost}. In particular, when $\Sigma=0$, we recover the deterministic PSD dynamics considered in \cite[Sec.~3.1]{Pates_Rantzer_2024}.

Furthermore, using the above lift, and with the condition $R_2\in\mathbb S_{>}^{n_u}$, the minimizer in \eqref{eq:generalminimizer_psd} evaluated at $X=xx^\top$ is of the form $U_1=xu^\top$ and $U_2=uu^\top$, with $u=a(x;\Lambda)=-Kx$, where $
K=M_{\theta,2}^{-1}(\Lambda)M_{\theta,1}^\top(\Lambda)=\left(B^\top\Lambda_{\theta,\Sigma}B+R_2\right)^{-1}\left(B^\top\Lambda_{\theta,\Sigma}A+R_1^\top\right)$.
Since $\Phi_\theta$ in \eqref{eq:Wishart_Phi} is independent of $p$, the LEQG problem admits the same recursion for $\Lambda$, with $\Gamma_\theta$ in \eqref{eq:Wishart_Gamma} evaluated at $p=1$. 

Note that, for $\theta>0$, writing $\Sigma=GG^\top$ and letting $\gamma^{-2}=2\theta$, the Riccati equation in \eqref{eq:Wishart_Phi} with the condition $I-2\theta\Lambda^{1/2}\Sigma\Lambda^{1/2}\in\mathbb S_{>}^{n_x}$ coincides with the Riccati equation associated with the corresponding $H_\infty$ problem (minimax problem) \cite{Glover_Doyle_1988,Basar_1991,Basar_Bernhard_2008}.

Other constructions of $\mrvW_t$ satisfying \eqref{eq:psd_w2condition} can be obtained by sampling affine processes on the PSD cone, including processes with jumps \cite{Cuchiero_etal_2011}, which have applications in optimal portfolio allocation under stochastic volatility and correlation risk \cite{Branger_etal_2017}.

\subsection{Wastewater Treatment with Risk Preferences}
\label{sec:wastewater_treatment}

Consider two wastewater holding tanks connected to a treatment reactor. Let $\mathcal X=\mathbb R_{\geq0}^3$, where $x^{(1)}$ and $x^{(2)}$ denote the pollutant mass in the tanks and $x^{(3)}$ denotes the pollutant mass in the reactor. The input $u\in\mathcal U(x):=[0,x^{(1)}]\times[0,x^{(2)}]$ specifies the mass transferred for treatment. The dynamics are $f(x,u,w)=Ax+Bu+[w^{(1)}~w^{(2)}~w^{(3)}]^\top$, where $A=\operatorname{diag}(4/5,4/5,1/5)$ and $B=\left[\begin{smallmatrix}-4/5&0\\0&-4/5\\1/5&1/5\end{smallmatrix}\right]$. Thus, $20\%$ of the remaining pollutant in each tank decays naturally per time step, while the reactor removes $80\%$ of its mass, including the transferred mass.

Let $\ell(x,u,w)=s^\top x+w^{(4)}$, where $s=[3/10~9/40~2/25]^\top$. The components of $\rvW\in\mathbb R_{\geq0}^4$ are conditionally independent, with $\rvW^{(i)}\sim\operatorname{Poisson}(\bar w^{(i)})$ for $i\in\{1,2,3\}$, where $\bar w\in\mathbb R_{>0}^3$, and $\rvW^{(4)}\mid(x,u)\sim\operatorname{Poisson}(u^{(1)}+u^{(2)})$. The first three components represent uncontrolled pollutant inflows, while $\rvW^{(4)}$ models costly filter replacements during treatment.

Using \eqref{eq:poisson_ex2_risk}, let $q_\theta:=(e^\theta-1)/\theta$ for $\theta\neq0$ and $q_0:=1$. Then $\mathcal D_\theta=\mathcal X^\ast=\mathbb{R}^{3}_{\geq0}$, $\Gamma_\theta(\lambda)=\bar w^\top\nu_{\theta,w}(\lambda)$, and $\Phi_\theta^{(i)}(\lambda)=s^{(i)}+\min\{\frac45\lambda^{(i)},\frac15\lambda^{(3)}+q_\theta\}$ for $i\in\{1,2\}$, while $\Phi_\theta^{(3)}(\lambda)=s^{(3)}+\frac15\lambda^{(3)}$. The unique fixed point satisfies $\lambda_\star^{(3)}=1/10$ and $\lambda_\star^{(i)}=\min\{5s^{(i)},s^{(i)}+1/50+q_\theta\}$ for $i\in\{1,2\}$. Hence, tank $i$ is treated if and only if $4s^{(i)}>1/50+q_\theta$. Since $\theta\mapsto q_\theta$ is continuous and strictly increasing with range $(0,\infty)$, let $\theta_i$ denote the unique solution of $q_{\theta_i}=4s^{(i)}-1/50$ for $i\in\{1,2\}$. The optimal stationary policy in \eqref{eq:opt_sol_star} is then
\begin{equation}
\label{eq:wastewater_policies}
a(x;\lambda_\star)=
\begin{cases}
0,&\theta>\theta_1\approx0.322,\\
[x^{(1)}~0]^\top,&\theta_2<\theta<\theta_1,\\
[x^{(1)}~x^{(2)}]^\top,&\theta<\theta_2\approx-0.261.
\end{cases}
\end{equation}
At $\theta=\theta_i$, every $a^{(i)}(x;\lambda_\star)\in[0,x^{(i)}]$ is optimal for tank $i$ for $i\in\{1,2\}$. Hence, sufficiently strong risk aversion favors natural decay to avoid uncertain treatment costs, risk neutrality treats only the first tank, and sufficiently strong risk seeking treats both tanks to retain the possibility of low treatment costs despite the risk of expensive filter replacements.

\section{Conclusion}
\label{sec:conclusion}
In this paper, we considered discrete-time stochastic optimal control problems with the state constrained to a proper cone while covering risk-neutral, risk-averse, and risk-seeking criteria within a single entropic-risk framework. We focused on problems whose Bellman operator maps linear functions to affine functions of the state. We characterized finite-horizon value functions and feedback policies through a recursion on the dual cone, and the optimal time-average cost and stationary policy through its smallest fixed point. Under finite-step cost observability, we established uniqueness, geometric convergence, and finite-horizon performance bounds for $\theta\leq0$. For $\theta\geq0$, we established uniform-in-time bounds on the expected state norm under stationary policies satisfying a Bellman dissipation inequality, with additional exponential-moment bounds for $\theta>0$. We also established an exact connection with deterministic infsup control for $\theta>0$ when the dynamics and costs are affine in the disturbances. Several examples illustrated the framework, including the connection between LQG and LQEG problems and a problem on the PSD cone. Future work will focus investigating conditions for uniqueness, convergence, and tractable convex formulations for $\theta>0$. Moreover, Further investigation of stability properties and moment bounds for both $\theta>0$ and $\theta<0$ remains an interesting direction for future work.

\appendices

\section{Proof of Lemma \ref{lem:finite_horizon_verification}}
\label{sec:finite_horizon_verification_proof}

\begin{proof}
For $t\in\{0,\ldots,N\}$, set
\begin{equation}
\label{eq:ver_value_functions}
V_t(x):=J_{N-t,\lambda}(x).
\end{equation}
Fix $\pi_N=\{a_t\}_{t=0}^{N-1}\in\Pi_N$. In what follows let $t\in\{0,\ldots,N-1\}$, write $\pi_N^t:=\{a_{t+s}\}_{s=0}^{N-t-1}$ and define
\begin{equation*}
V_t^{\pi_N}(x):=
\mathcal R_\theta\left(S_{N-t,\pi_N^t}^\lambda(x)\right),\quad
V_N^{\pi_N}(x):=\lambda^\top x.
\end{equation*}
Then
\begin{align}
\label{eq:ver_endpoints}
V_0^{\pi_N}(x)&=J_{N,\lambda}^{\pi_N}(x),\\
\label{eq:ver_terminal_values}
V_N^{\pi_N}(x)&=\lambda^\top x=V_N(x).
\end{align}

By Assumption~\ref{ass:basic}, admissibility, and
$\lambda\in\mathcal X^\ast$, we have
$0\leq S_{N-t,\pi_N^t}^\lambda(x)<\infty$
almost surely. Hence its expectation and exponential
moments are well defined in the extended sense, in particular, for $\theta\neq0$
\begin{equation*}
\mathbb E\left[e^{\theta S_{N-t,\pi_N^t}^\lambda(x)}\right]
\in\begin{cases}
[1,+\infty],&\theta>0,\\
(0,1],&\theta<0.
\end{cases}
\end{equation*}
and
\begin{align}\label{etoE}
e^{\theta V_t^{\pi_N}(x)}
=\mathbb E\left[
e^{\theta S_{N-t,\pi_N^t}^\lambda(x)}\right].
\end{align}
Now using \eqref{eq:disturbance}, the expectation can be written as
\begin{multline}
\label{eq:ver_iterated_exponential}
\mathbb E\left[
e^{\theta S_{N-t,\pi_N^t}^\lambda(x)}\right]\\
=\int_{\mathcal W}\cdots\int_{\mathcal W}
\exp\Bigl(\theta\Bigl[
\sum_{k=t}^{N-1}\ell(x_k,u_k,w_k)
+\lambda^\top x_N\Bigr]\Bigr)\\
\,\rcd(dw_{N-1}\mid x_{N-1},u_{N-1})
\cdots \rcd(dw\mid x,u),
\end{multline}
with $w_t=w$, $x_t=x$ and $u=a_t(x)$. Separating $\ell(x,u,w)$ in
\eqref{eq:ver_iterated_exponential}, the remaining
integrals give $e^{\theta V_{t+1}^{\pi_N}(f(x,u,w))}$ by \eqref{etoE}. Therefore,
\begin{equation}
\label{eq:ver_exponential_recursion}
e^{\theta V_t^{\pi_N}(x)}
=\int_{\mathcal W}e^{\theta\ell(x,u,w)}
\, e^{\theta V_{t+1}^{\pi_N}(f(x,u,w))}
\,\rcd(dw\mid x,u).
\end{equation}
For $\theta=0$, the same successive integration,
applied to the accumulated cost, gives
\begin{equation}
\label{eq:ver_expectation_recursion}
V_t^{\pi_N}(x)
=\int_{\mathcal W}\Bigl(\ell(x,u,w)
+V_{t+1}^{\pi_N}(f(x,u,w))\Bigr)
\,\rcd(dw\mid x,u).
\end{equation}
Taking logarithms and dividing by $\theta\neq0$ in
\eqref{eq:ver_exponential_recursion}, or using
\eqref{eq:ver_expectation_recursion} for $\theta=0$,
we obtain
\begin{equation}
\label{eq:ver_policy_recursion}
V_t^{\pi_N}(x)
=\mathcal R_\theta\Bigl(\ell(x,u,\rvW)
+V_{t+1}^{\pi_N}(f(x,u,\rvW))\mid x,u\Bigr).
\end{equation}

We now prove $V_t^{\pi_N}(x)\geq V_t(x)$ by backward
induction. For $t=N$, equality holds by
\eqref{eq:ver_terminal_values}.
Suppose $V_{t+1}^{\pi_N}(x)\geq V_{t+1}(x)$ for every
$x\in\mathcal X$. By \eqref{eq:ver_policy_recursion},
admissibility, \eqref{eq:monoticity_of_R}, and
\eqref{eq:rs_JN1},
\begin{multline}
\label{eq:ver_cost_comparison}
V_t^{\pi_N}(x)\geq\mathcal R_\theta\Bigl(\ell(x,u,\rvW)
+V_{t+1}(f(x,u,\rvW))\mid x,u\Bigr)\\
\geq(T_\theta V_{t+1})(x)=V_t(x).
\end{multline}
Thus, by induction, \eqref{eq:ver_value_functions} and \eqref{eq:ver_endpoints},
\begin{equation}
\label{eq:ver_policy_lower_bound}
J_{N,\lambda}^{\pi_N}(x)\geq J_{N,\lambda}(x),~x\in\mathcal X.
\end{equation}

Now let $\pi_N^\star=\{a_t^\star\}_{t=0}^{N-1}$
satisfy \eqref{eq:rs_optimal_feedback}.
Since its decision rules are Borel measurable and
admissible, $\pi_N^\star\in\Pi_N$.
Starting from $V_N^{\pi_N^\star}=V_N$, the same
induction, together with \eqref{eq:rs_optimal_feedback},
gives
\begin{multline*}
V_t^{\pi_N^\star}(x)=\mathcal R_\theta\Bigl(\ell(x,a_t^\star(x),\rvW)\\
+V_{t+1}(f(x,a_t^\star(x),\rvW))
\mid x,a_t^\star(x)\Bigr)
=(T_\theta V_{t+1})(x)=V_t(x).
\end{multline*}
Hence $J_{N,\lambda}^{\pi_N^\star}(x)=J_{N,\lambda}(x)$.
Combining this with \eqref{eq:ver_policy_lower_bound}
gives
\begin{align*}
J_{N,\lambda}(x)
=\inf_{\pi_N\in\Pi_N}J_{N,\lambda}^{\pi_N}(x)
=J_{N,\lambda}^{\pi_N^\star}(x),
~ x\in\mathcal X,
\end{align*}
which proves \eqref{eq:rs_JN2}.
\end{proof}

\bibliographystyle{IEEEtran} 
\bibliography{IEEEabrv,litabrv} 

\begin{thebibliography}{10}
\providecommand{\url}[1]{#1}
\csname url@samestyle\endcsname
\providecommand{\newblock}{\relax}
\providecommand{\bibinfo}[2]{#2}
\providecommand{\BIBentrySTDinterwordspacing}{\spaceskip=0pt\relax}
\providecommand{\BIBentryALTinterwordstretchfactor}{4}
\providecommand{\BIBentryALTinterwordspacing}{\spaceskip=\fontdimen2\font plus
\BIBentryALTinterwordstretchfactor\fontdimen3\font minus \fontdimen4\font\relax}
\providecommand{\BIBforeignlanguage}[2]{{%
\expandafter\ifx\csname l@#1\endcsname\relax
\typeout{** WARNING: IEEEtran.bst: No hyphenation pattern has been}%
\typeout{** loaded for the language `#1'. Using the pattern for}%
\typeout{** the default language instead.}%
\else
\language=\csname l@#1\endcsname
\fi
#2}}
\providecommand{\BIBdecl}{\relax}
\BIBdecl

\bibitem{Colaneri_Middleton_Blanchini_2016}
P.~Colaneri, R.~H. Middleton, and F.~Blanchini, ``Optimal control of a class of positive {Markovian} bilinear systems,'' \emph{Nonlinear Anal. Hybrid Syst.}, vol.~21, pp. 155--170, Aug. 2016.

\bibitem{Dhingra_Colombino_Jovanovic_Rantzer_Smith_2016}
\BIBentryALTinterwordspacing
N.~K. Dhingra, M.~Colombino, M.~R. Jovanovi{\'c}, A.~Rantzer, and R.~S. Smith, ``On the optimal control problem for a class of monotone bilinear systems,'' 2016. [Online]. Available: \url{https://arxiv.org/abs/1611.09984}
\BIBentrySTDinterwordspacing

\bibitem{Zorzan_Rantzer_2017}
I.~Zorzan and A.~Rantzer, ``{$L_1$} and {$H_\infty$} optimal control of positive bilinear systems,'' in \emph{Proc. 56th Annu. IEEE Conf. Decis. Control (CDC)}.\hskip 1em plus 0.5em minus 0.4em\relax IEEE, 2017, pp. 727--732.

\bibitem{Rantzer_2022}
A.~Rantzer, ``Explicit solution to {Bellman} equation for positive systems with linear cost,'' in \emph{Proc. 61st IEEE Conf. Decis. Control (CDC)}.\hskip 1em plus 0.5em minus 0.4em\relax IEEE, 2022, pp. 6154--6155.

\bibitem{DynamicProgramming_YuchaoRantzer}
Y.~Li and A.~Rantzer, ``Exact dynamic programming for positive systems with linear optimal cost,'' \emph{{IEEE} Trans. Autom. Control}, vol.~69, no.~12, pp. 8738--8750, 2024.

\bibitem{Ohlin_Pates_Arcak_2025}
D.~Ohlin, R.~Pates, and M.~Arcak, ``On exact solutions to the linear {Bellman} equation,'' \emph{IEEE Control Syst. Lett.}, vol.~9, pp. 1568--1573, 2025.

\bibitem{semilineardynamicprogramminganalysis}
Y.~Li and D.~P. Bertsekas, ``Semilinear dynamic programming: Analysis, algorithms, and certainty equivalence properties,'' \emph{{IEEE} Trans. Autom. Control}, 2026, in press.

\bibitem{Pates_Rantzer_2024}
R.~Pates and A.~Rantzer, ``Optimal control on positive cones,'' in \emph{Proc. 63rd IEEE Conf. Decis. Control (CDC)}.\hskip 1em plus 0.5em minus 0.4em\relax IEEE, 2024, pp. 992--997.

\bibitem{Glover_Doyle_1988}
K.~Glover and J.~C. Doyle, ``State-space formulae for all stabilizing controllers that satisfy an {$H_\infty$}-norm bound and relations to risk sensitivity,'' \emph{Syst. Control Lett.}, vol.~11, no.~3, pp. 167--172, 1988.

\bibitem{Basar_Bernhard_2008}
T.~Ba{\c{s}}ar and P.~Bernhard, \emph{{$H_\infty$}-Optimal Control and Related Minimax Design Problems: A Dynamic Game Approach}, 2nd~ed., ser. Modern Birkh{\"a}user Classics.\hskip 1em plus 0.5em minus 0.4em\relax Boston, MA: Birkh{\"a}user, 2008.

\bibitem{Gurpegui_Tegling_Rantzer_2023}
A.~Gurpegui, E.~Tegling, and A.~Rantzer, ``Minimax linear optimal control of positive systems,'' \emph{IEEE Control Syst. Lett.}, vol.~7, pp. 3920--3925, 2023.

\bibitem{Gurpegui_Jeeninga_Tegling_Rantzer_2026}
A.~Gurpegui, M.~Jeeninga, E.~Tegling, and A.~Rantzer, ``Minimax linear regulator problems for positive systems,'' \emph{{IEEE} Trans. Autom. Control}, vol.~71, no.~8, pp. 5290--5303, Aug. 2026.

\bibitem{Rantzer_2025_AdaptiveCones}
A.~Rantzer, ``Minimax optimal adaptive control for systems on cones,'' in \emph{Proc. 64th IEEE Conf. Decis. Control (CDC)}, 2025, pp. 4137--4139.

\bibitem{Whittle_2002}
P.~Whittle, ``Risk sensitivity, a strangely pervasive concept,'' \emph{Macroecon. Dyn.}, vol.~6, no.~1, pp. 5--18, 2002.

\bibitem{SDEAlba}
A.~Gurpegui, T.~Tanaka, and A.~Rantzer, ``{$L_1$} optimal control of continuous-time stochastic positive systems,'' \emph{IEEE Control Syst. Lett.}, vol.~10, pp. 1573--1578, 2026.

\bibitem{Howard_Matheson_1972}
R.~A. Howard and J.~E. Matheson, ``Risk-sensitive {Markov} decision processes,'' \emph{Manage. Sci.}, vol.~18, no.~7, pp. 356--369, 1972.

\bibitem{Whittle_1990}
P.~Whittle, \emph{Risk-Sensitive Optimal Control}, ser. Wiley-Interscience Series in Systems and Optimization.\hskip 1em plus 0.5em minus 0.4em\relax Chichester: Wiley, 1990, vol.~2.

\bibitem{Jacobson_1973}
D.~H. Jacobson, ``Optimal stochastic linear systems with exponential performance criteria and their relation to deterministic differential games,'' \emph{{IEEE} Trans. Autom. Control}, vol.~18, no.~2, pp. 124--131, 1973.

\bibitem{Whittle_1981}
P.~Whittle, ``Risk-sensitive linear/quadratic/{G}aussian control,'' \emph{Adv. Appl. Probab.}, vol.~13, no.~4, pp. 764--777, 1981.

\bibitem{Bauerle_Rieder_2014}
N.~B{\"a}uerle and U.~Rieder, ``More risk-sensitive {Markov} decision processes,'' \emph{Math. Oper. Res.}, vol.~39, no.~1, pp. 105--120, Feb. 2014.

\bibitem{Chapman_RADP2022}
M.~P. Chapman and K.~M. Smith, ``Classical risk-averse control for a finite-horizon {Borel} model,'' \emph{IEEE Control Syst. Lett.}, vol.~6, pp. 1525--1530, 2022.

\bibitem{Jaskiewicz_2007}
A.~Ja{\'s}kiewicz, ``Average optimality for risk-sensitive control with general state space,'' \emph{Ann. Appl. Probab.}, vol.~17, no.~2, pp. 654--675, 2007.

\bibitem{Campi_James_1996}
M.~C. Campi and M.~R. James, ``Nonlinear discrete-time risk-sensitive optimal control,'' \emph{Int. J. Robust Nonlinear Control}, vol.~6, no.~1, pp. 1--19, 1996.

\bibitem{Dupuis_James_Petersen_2000}
P.~Dupuis, M.~R. James, and I.~R. Petersen, ``Robust properties of risk-sensitive control,'' \emph{Math. Control Signals Syst.}, vol.~13, no.~4, pp. 318--332, 2000.

\bibitem{Hansen_Sargent_2008}
L.~P. Hansen and T.~J. Sargent, \emph{Robustness}.\hskip 1em plus 0.5em minus 0.4em\relax Princeton, NJ: Princeton University Press, 2008.

\bibitem{BoydVandenberghe2004}
S.~Boyd and L.~Vandenberghe, \emph{Convex Optimization}.\hskip 1em plus 0.5em minus 0.4em\relax Cambridge University Press, 2004.

\bibitem{Bertsekas_Shreve_1996}
D.~P. Bertsekas and S.~E. Shreve, \emph{Stochastic Optimal Control: The Discrete-Time Case}, ser. Optimization and Neural Computation Series.\hskip 1em plus 0.5em minus 0.4em\relax Belmont, MA: Athena Scientific, 1996, vol.~5.

\bibitem{Goldstein_2011}
D.~H. Goldstein, \emph{Polarized Light}, 3rd~ed.\hskip 1em plus 0.5em minus 0.4em\relax Boca Raton, FL: CRC Press, 2011.

\bibitem{Boxma_Ivanovs_Kosinski_Mandjes_2011}
O.~Boxma, J.~Ivanovs, K.~Kosi{\'n}ski, and M.~Mandjes, ``L{\'e}vy-driven polling systems and continuous-state branching processes,'' \emph{Stoch. Syst.}, vol.~1, no.~2, pp. 411--436, 2011.

\bibitem{Rantzer_2006}
A.~Rantzer, ``Relaxed dynamic programming in switching systems,'' \emph{IEE Proc. Control Theory Appl.}, vol. 153, no.~5, pp. 567--574, 2006.

\bibitem{Piotrowski_Cavalcante_2022}
T.~Piotrowski and R.~L.~G. Cavalcante, ``The fixed point iteration of positive concave mappings converges geometrically if a fixed point exists: Implications to wireless systems,'' \emph{{IEEE} Trans. Signal Process.}, vol.~70, pp. 4697--4710, 2022.

\bibitem{Bertsekas_2022_ADP}
D.~P. Bertsekas, \emph{Abstract Dynamic Programming}, 3rd~ed.\hskip 1em plus 0.5em minus 0.4em\relax Athena Scientific, 2022.

\bibitem{Vysotsky_2021}
V.~Vysotsky, ``When is the rate function of a random vector strictly convex?'' \emph{Electron. Commun. Probab.}, vol.~26, pp. 1--11, 2021.

\bibitem{Gurpegui_Tegling_Rantzer_2025_Unbounded}
\BIBentryALTinterwordspacing
A.~Gurpegui, E.~Tegling, and A.~Rantzer, ``A minimax optimal controller for positive systems,'' 2025. [Online]. Available: \url{https://arxiv.org/abs/2502.01180}
\BIBentrySTDinterwordspacing

\bibitem{Anderson_1946}
T.~W. Anderson, ``The non-central {W}ishart distribution and certain problems of multivariate statistics,'' \emph{Ann. Math. Statist.}, vol.~17, no.~4, pp. 409--431, 1946.

\bibitem{Muirhead_1982}
R.~J. Muirhead, \emph{Aspects of Multivariate Statistical Theory}.\hskip 1em plus 0.5em minus 0.4em\relax New York: John Wiley \& Sons, 1982.

\bibitem{Basar_1991}
T.~Ba{\c{s}}ar, ``A dynamic games approach to controller design: Disturbance rejection in discrete-time,'' \emph{{IEEE} Trans. Autom. Control}, vol.~36, no.~8, pp. 936--952, 1991.

\bibitem{Cuchiero_etal_2011}
C.~Cuchiero, D.~Filipovi{\'c}, E.~Mayerhofer, and J.~Teichmann, ``Affine processes on positive semidefinite matrices,'' \emph{Ann. Appl. Probab.}, vol.~21, no.~2, pp. 397--463, 2011.

\bibitem{Branger_etal_2017}
N.~Branger, M.~Muck, F.~T. Seifried, and S.~Weisheit, ``Optimal portfolios when variances and covariances can jump,'' \emph{J. Econ. Dyn. Control}, vol.~85, pp. 59--89, 2017.

\end{thebibliography}
\vfill
\end{document}